\documentclass[11pt,reqno]{article}
\usepackage{eurosym}
\usepackage{comment}
\usepackage{mathptmx}
\usepackage{amssymb}
\usepackage{amsfonts}
\usepackage{amsthm}
\usepackage{amsmath}
\usepackage{dsfont}
\usepackage{graphicx}
\usepackage{shadow}
\usepackage[all]{xy}
\usepackage{enumerate}
\usepackage{makeidx}
\usepackage{mathrsfs}
\usepackage{upgreek}
\usepackage{stmaryrd}
\usepackage[utf8]{inputenc}
\usepackage{array}
\usepackage{tikz-cd}
\usepackage{thumbpdf}
\usepackage[pagebackref]{hyperref}
\usepackage{textcomp}
\usepackage{chapterbib}
\usepackage{graphics}
\usepackage{longtable}
\usepackage{epsf}
\usepackage{bm}
\usepackage{adjustbox}
\usepackage{centernot}
\usepackage{setspace}
\usepackage{imakeidx}
\usepackage{tikz}
\usepackage{tikz-cd}
\usepackage{float}
\usepackage{caption}

\makeindex
\newtheorem{theorem}{Theorem}[section]
\newtheorem{lem}[theorem]{Lemma}
\newtheorem{thm}[theorem]{Theorem}
\newtheorem{prop}[theorem]{Proposition}
\newtheorem{cor}[theorem]{Corollary}
\theoremstyle{definition}
\newtheorem*{notatation}{Notation}
\newtheorem*{Beweis}{Proof}
\newtheorem{defn}[theorem]{Definition}

\newtheorem{defns}[theorem]{Definitions}
\newtheorem{rem}[theorem]{Remark}

\newtheorem{punto}[theorem]{}
\newtheorem*{notation}{Notation}

\theoremstyle{remark}

\newtheorem{ex}[theorem]{Example}
\newtheorem{exs}[theorem]{Examples}
\newtheorem{remarks}[theorem]{Remarks}
\begin{document}

\title{Higher Separation Axioms for $X$-top Lattices\\
Applications to Commutative (Semi)rings\thanks{%
MSC2020: Primary: 06B23 , 6B30, 6B35; Secondary: 13C13, 16Y60} \thanks{%
\textbf{Keywords: }Lattices, $X$-Top Lattice, Zariski Topology, Separation
Axioms, Regular Spaces, Completely Regular Spaces, Normal Spaces, Completely
Normal Spaces, Perfectly Normal Spaces, Zero Dimensional (Semi)Rings}
{\small \thanks{\text{Extracted from the M.Sc. thesis of the second author
under the supervision of Prof. J. Abuhlail.}}}}
\author{$%
\begin{array}{cc}
\text{Jawad Abuhlail}{\small \thanks{\text{Corresponding Author}}} & \text{%
Abdulmushin Alfaraj} \\
\text{{\small abuhlail@kfupm.edu.sa}} & \text{{\small %
abdulmuhsinalfaraj@hotmail.com}} \\
\text{{\small Department of Mathematics}} & \text{{\small Department of
Mathematical Sciences}} \\
\text{{\small King Fahd University of Petroleum $\&$ Minerals}} & \text{%
{\small University of Bath}} \\
\text{{\small 31261 Dhahran, KSA}} & \text{{\small Celeverton Down, Bath,
BA27AY, UK}}%
\end{array}%
$}
\date{\today }
\maketitle

\begin{abstract}
We study several separation axioms for $X$-top-lattices (i.e. a lattice $L$
for which a given subset $X\subseteq L\backslash \{1\}$ admits a \emph{%
Zariski-like topology}). Such spaces are $T_{0}$ and usually far away from
being $T_{2}.$ We provide sufficient/necessary conditions for an $X$-top
lattice so that $X$ is $T_{2},$ \emph{regular} ($T_{3}$), \emph{completely
regula}r ($T_{3\frac{1}{2}}$), \emph{normal}, \emph{completely normal} or
\emph{perfectly normal} ($T_{6}$). We apply our results mainly to the
spectrum of prime (resp. maximal, minimal) ideals of a commutative
(semi)ring. We illustrate our results with several examples/counterexamples.
\end{abstract}

\section*{Introduction}

Different spectra of \emph{prime submodules }associated to a module $M$ over
a ring $R$ were investigated (cf. \cite{MMS1997}, \cite{Wij2006}). Several
\emph{Zariski-like topologies} were associated to these spectra (e.g., \cite%
{Lu1999}, \cite{MS2006}) and investigated by many authors (e.g., \cite%
{BH2008}, \cite{Tek2009}) including the first author (e.g., \cite{Abu2011-CA}%
). These were also dualized to what we call \emph{dual Zariski-like
topologies} on the spectrum of \emph{coprime submodules} of a module $M$
over an associative ring $\mathcal{A}$ (e.g., \cite{Abu2011}, \cite{Abu2015}%
) as well as the spectrum of \emph{coprime subcomodules} of a comodule $N$
over a coring $\mathcal{C}$ (e.g., \cite{Abu2008}, \cite{Abu2006}).
Moreover, such investigations were carried out to the spectrum of
prime/coprime submodules of a semimodule over a semiring (e.g., \cite%
{AUT2013}, \cite{HPH2021}).

\bigskip

The $X$\emph{-top lattices} were introduced by Abuhlail and Lomp \cite%
{AL2016} as a general framework for investigating (dual) Zariski-like
topologies on the spectra of (co)prime modules/comodules over (co)rings as
well as the spectra of (co)prime (semi)modules over (semi)rings. An
intensive study of the properties of such topologies is \cite{AH2019}, where
the main interest was investigating the \emph{interplay} between the
algebraic properties of the algebraic structures considered and the
topological properties of the (dual) Zariski-like topologies associated to
them.

\bigskip

This paper is a continuation of \cite{AA} in which we studied the \emph{%
lower separation axioms }$T_{\frac{1}{4}},$ $T_{\frac{1}{2}},$ $T_{\frac{3}{4%
}}$ and $T_{1}$ for $X$-top lattices. In this paper, we study the $X$-top
lattices for which $X$ is $T_{2}$ as well as those for which $X$ is \emph{%
regular}, \emph{completely regular}, \emph{normal}, \emph{completely normal}
or \emph{perfectly normal}. Following \cite{Wil1970}, we use the versions of
these regularity/normality properties that do \emph{not} assume $T_{1}$
since this separation axiom forces $X$ to be \emph{Krull zero-dimensional}
\cite[Proposition 2.9.]{AA}. However, for each $X$-top lattice $L,$ we
have:\ $X$ is $T_{0},$ whence $X$ is regular (resp. completely regular,
perfectly normal) if and only if $X$ is $T_{3}$ (resp. $T_{3\frac{1}{2}},$ $%
T_{6}$).

\bigskip

In Section 1, we recall some preliminaries from the theory of lattices \cite%
{Gra2010}, $X$-top lattices \cite{AL2016}, the theory of semirings \cite%
{Gol1999}, as well as general topology \cite{Wil1970}).

\bigskip

In Section 2, we study $X$-top lattices for which $X$ is normal, completely
normal and perfectly normal. We investigate, in particular, the relation
between the \textbf{normality} of $X,$ the \textbf{pm-property} of $X$
(i.e., every element of $X$ is comparable with a \emph{unique} $\mathfrak{m}%
\in Max(X)$) and the \textbf{max-retractibility} of $X$ (i.e., the existence
of a \emph{retraction} $\mu :Max(X)\longrightarrow X$). While these three
properties are equivalent for spectral spaces (cf. Proposition \ref{tspm}),
we show in Example \ref{pm-not-normal} that $Y=Max(Spec(\mathbb{Z})),$ which
has \emph{trivially} the pm-property and the max-retractibility, is far away
from being normal as it is \emph{extremely non-normal} (even \emph{extremely
non-Hausdorff}). In Proposition \ref{txnpm}, we show that in case $X$ is
coatomic: the normality of $X$, as well as the max-retractability of $X,$ is
a sufficient condition for $X$ to have the pm-property. The converse is
investigated in Theorem \ref{pxn}. In Theorem \ref{FMax-pm-retract}, we show
that \emph{all} these three properties are equivalent in case $X$ is
coatomic, atomic with both $Min(X)\ $and $Max(X)$ finite. This is in,
particular, the case when $X$ is finite (cf. Corollary \ref%
{finite-pm-retract}).

Then we investigate the interplay between the locality of $X$ (i.e., the
existence of a \emph{unique} $\mathfrak{m}\in X$ such that $x\leq \mathfrak{m%
}$ for all $x\in X$) and its normality. We show in Theorem \ref{cor-local}
that, in case $X$ is coatomic and colocal (e.g., $\sqrt[X]{0_{L}}\in X$),
the locality of $X$ is equivalent to the normality of $X,$ which is in turn
equivalent to the hyperconnectedness of $X.$ This applies, in particular, to
$X=Spec(R),$ where $R$ is an integral (semi)domain (cf. Corollary \ref%
{domain-normal}).

Theorem \ref{tcn} provides a sufficient/necessary condition for such $X$ to
be \emph{completely normal} in case $X$ is a forest consisting of a finite
number of \emph{strongly disjoint} $\bigwedge $-trees with finite base or
contains a $\bigvee $-tree $\mathcal{V}$ with a finite cover. Example \ref%
{strongly-disjoint} illustrates this result while Example \ref{not-pm}
demonstrates that the assumption that the $\bigwedge $-trees in Theorem \ref%
{tcn} (1) are strongly disjoint cannot be weakened by replacing it with the
\emph{standard disjointness} of the $\bigwedge $-trees involved as (in this
case, $X$ might lack the normality and not only then complete normality).

In Section 3, we study $X$-top lattices for which $X$ is \emph{regular}
(equivalently $T_{3}$). In Proposition \ref{cpmpact-T3-T4}, we show in case $%
X$ is \emph{compact}: the regularity of $X$ is equivalent to $T_{4}$ (and $%
T_{2}$). In Theorem \ref{abs-flat-regular}, we characterize the $X$-top
lattices for which $X$ is a \emph{Stone space} as those for which $X$ is
spectral and regular (equivalently, $X$ is homeomorphic to $Spec(R)$ for
some \emph{Jacobson} \emph{pm-(semi)ring}). As a consequence, Corollary \ref%
{vN-regular} characterizes the commutative von Neumann regular rings as the
reduced ones with regular prime spectra, or equivalently as the reduced\emph{%
\ }(\emph{dual})\emph{\ Jacobson} commutative rings with normal prime
spectra. In Example \ref{B(3,1)}, we provide an example demonstrating that
these characterizations are not valid for \emph{proper} semirings (that are
\emph{not} rings).

Motivated by the breaking result of W. Lewis \cite{Lew1973}, which states
roughly that every \emph{finite} poset can be realized as the prime spectrum
of a some commutative ring, we devote the last part of this paper to
studying the regularity and normality of several \emph{finite} posets for
which a given lattice $\mathcal{L}$ is a \emph{top lattice}. Moreover, we
realize some of these posets as the prime spectra of some \emph{proper}
semirings. Moreover, we demonstrate that the \emph{regularity} in this
context is stronger than the \emph{normality} in the following sense:\ while
the regularity of compact $X$ implies the normality of $X$ (by Proposition %
\ref{cpmpact-T3-T4}) we have examples of $X$-top lattices for which $X$ is
finite (whence spectral) and normal but \emph{not} regular. In Example \ref%
{normal-not-regular}, $X$ is normal but \emph{extremely non-regular}. In
Example \ref{CN-not-PN}, $Y$ is even $T_{\frac{3}{4}}$ and completely normal
but \emph{extremely non-regular}. In Example \ref{anti-normal}, each of
these examples, $X$ and $Y$ are normal but \emph{not} perfectly normal as it
contains a $\mathcal{C}_{2}$ (cf. Proposition \ref{ppn}).

\section{Preliminaries}

\subsection*{Lattices}

We recall some relevant definitions and results from Lattice Theory. We
follow \cite{Gra2010} (unless otherwise stated explicitly).

\begin{punto}
Let $(L,\wedge )$ be a \emph{complete} meet-semilattice and $B,C\subseteq L.$
We say that $q\in C$ is \textbf{(completely) strongly }$B$-\textbf{%
irreducible in }$L$ iff for any $A\underset{\text{finite}}{\subseteq }B$ ($%
A\subseteq B$), we have: $\bigwedge\limits_{a\in A}a\leq q\Longrightarrow
a\leq q$ for some $a\in A.$

With $SI^{B}(C)$ (resp. $CSI^{B}(C)$), we denote the set of strongly $B$%
-irreducible (resp. completely strongly $B$-irreducible) elements of $C.$ We
drop the superscript $B$ if it is clear from the context.
\end{punto}

\begin{punto}
Let $\mathcal{L}=(L;\vee ,0;\wedge ,1)$ be a \emph{bounded} lattice and $%
\emptyset \neq X\subseteq A\subseteq L.$ For $a\in L,$ we set%
\begin{equation*}
Max(a;X):=\{\mathfrak{m}\in Max(X)\mid a\leq \mathfrak{m}\}\text{ and }%
Min(a;X):=\{m\in Min(X)\mid m\leq a\}.
\end{equation*}%
We say that $A$ is

$X$-\textbf{atomic} iff for every $a\in A:$ there exists $m\in Min(X)$ such
that $m\leq a;$

$X$-\textbf{coatomic} iff for every $a\in A:$ there exists $\mathfrak{m}\in
Max(X)$ such that $a\leq \mathfrak{m}.$

We say that $\emptyset \neq X\subseteq L$ is

\textbf{atomic} (\textbf{coatomic}) iff $X$ is $X$-atomic ($X$-coatomic);

\textbf{coocal} (\textbf{local}) iff $X$ is atomic and $\left\vert
Min(X)\right\vert =1$ (coatomic and $Max(X)=1$).
\end{punto}

\begin{remarks}
Let $\mathcal{L}=(L;\vee ,0;\wedge ,1)$ be a \emph{bounded} lattice.

\begin{enumerate}
\item If $\emptyset \neq X\subseteq L$ is finite, then $X$ is atomic and
coatomic.

\item Our notion of locality for $X$ should not be confused with that of a
\emph{local lattice} (a \textbf{frame}), i.e., a complete lattice in which
meets distribute over arbitrary joins (cf. \cite[page 7]{Gol1999}).
\end{enumerate}
\end{remarks}

\subsection*{$X$-top Lattices}


We recall some definitions and notation from the \emph{Theory }$X$\emph{-Top
Lattices.} We follow \cite{AL2016} (unless otherwise stated explicitly).

\begin{notation}
Let $\mathcal{L}=(L;\vee ,0;\wedge ,1)$ be a complete lattice and $\emptyset
\neq X\subseteq L\backslash \{1\}$. For any $a\in L$, define%
\begin{equation*}
V_{X}(a):=\{x\in X\ |\ a\leq x\},\text{ }D_{X}(a):=X\backslash V_{X}(a)\text{
and }U_{X}(a):=\{x\in X\ |\ x\leq a\}.
\end{equation*}%
We call $V_{X}(a)$ the \emph{variety} of $a$ in $X.$ Moreover, we set
\begin{equation*}
V_{X}(\mathcal{L}):=\{V_{X}(a)\ |\ a\in L\}\text{ and }\tau _{X}(\mathcal{L}%
)=\{D_{X}(a)\ |\ a\in L\}.
\end{equation*}%
We drop the subscript $X$ if it is clear from the context.
\end{notation}

\begin{defn}
Let $\mathcal{L}=(L;\vee ,0;\wedge ,1)$ be a complete lattice and $\emptyset
\neq X\subseteq L\backslash \{1\}.$ We say that $\emptyset \neq A,B\subseteq
X$ are strongly disjoint iff $V(\bigwedge\limits_{a\in A}a)\cap
V(\bigwedge\limits_{b\in B}b)=\emptyset .$
\end{defn}

\begin{punto}
\label{IV-def}Let $\mathcal{L}=(L;\vee ,0;\wedge ,1)$ be a \emph{complete
lattice} and $\emptyset \neq X\subseteq L\backslash \{1\}.$ Notice that $%
V(0)=X,$ $V(1)=\emptyset $ and $V(\mathcal{L})$ is closed under arbitrary
intersections as $\bigcap_{a\in A}(V(a))=V(\bigvee_{a\in A}a)$ for any $%
A\subseteq L.$ We say that $\mathcal{L}$ is an $X$\textbf{-top lattice} \cite%
{AL2016} iff $V(\mathcal{L})$ is closed under finite unions.\newline
Consider $\emptyset \neq X\subseteq L\backslash \{1\}$. For any $Y\subseteq
X $ and $a\in L,$ we set%
\begin{equation*}
I_{X}(Y):=\bigwedge_{y\in Y}y\text{ and }\sqrt[X]{a}:=I_{X}(V_{X}(a))\text{
and }C^{X}(\mathcal{L}):=\{a\in L\text{ }|\text{ }a=\sqrt[X]{a}\}.
\end{equation*}%
We drop the superscript $X$ if it is clear from the context. Clearly, $%
\emptyset \neq X\subseteq C^{X}(\mathcal{L})$ and $(C^{X}(\mathcal{L}%
),\wedge )$ is a meet-semilattice.
\end{punto}

Now, we recall a fundamental characterization of $X$-top lattices by
Abuhlail and Lomp \cite{AL2016}.

\begin{thm}
\label{xct}\emph{(\cite[\emph{Theorem 2.2}]{AL2016})} Let $\mathcal{L}%
=(L;\vee ,0;\wedge ,1)$ be a complete lattice and $\emptyset \neq X\subseteq
L\backslash \{1\}$. Then $\mathcal{L}$ is an $X$-top lattice if and only if $%
X=SI^{C^{X}(\mathcal{L})}(X).$
\end{thm}

The following is a direct, but very useful, consequence of Theorem \ref{xct}
especially in constructing examples and counterexamples.

\begin{cor}
\label{Y-top}\emph{(\cite[Corollary 1.9]{AA})} Let $\mathcal{L}=(L;\vee
,0;\wedge ,1)$ be an $X$-top lattice for some $\emptyset \neq X\subseteq
L\backslash \{1\}.$ If $\emptyset \neq Y\subseteq X,$ then $L$ is a $Y$-top
lattice and the corresponding topology on $Y$ is the induced subspace
topology.
\end{cor}

\subsection*{Semirings}

We recall some definitions and examples from the Theory of Semirings. We
follow \cite{Gol1999} (unless otherwise stated explicitly).

\begin{punto}
A \textbf{semiring} $R$ is roughly a ring \emph{not necessarily} with
subtraction. We assume hat $0_{R}$ is \emph{absorbing} (i.e., $0\cdot
r=0=r\cdot 0$ for all $r\in R$) and that $0_{R}\neq 1_{R}.$ If, in addition,
the monoid $(R,\cdot )$ is commutative, we say that $R$ is a \emph{%
commutative semiring}. We call a (commutative) semiring with no non-zero
zerodivisors \textbf{entire} \cite{Gol1999} (\textbf{semidomain}). An entire
semiring (semidomain) in which every non-zero element has a multiplicative
inverse is called a \textbf{division semiring} (\textbf{semifield}). A
semiring that is not a ring is called a \textbf{proper semiring.}
\end{punto}

\begin{ex}
\label{B(n,i)}(\cite{AA1994}, \cite[Example 1.8]{Gol1999}) Consider $%
B(n,i):=({\{0,1,2,...,n-1\}},\oplus ,0,\odot ,1),$ where

\begin{enumerate}
\item $x\oplus y=x+y\ $ if $\ x+y\leq n-1$; otherwise, $x+y=u,$ the unique
positive integer satisfying $i\leq u\leq n-1$ and $x+y\equiv u\ $mod\ $(n-i)$%
;

\item $x\odot y=xy\ $ if $\ xy\leq n-1$; otherwise, $xy=v$ the unique
positive integer satisfying $i\leq v\leq n-1$ and $xy\equiv v\ $mod\ $(n-i)$.
\end{enumerate}

Then $B(n,i)$ is a commutative semiring. Observe that $B(2,1)=\mathbb{B}$,
the \textbf{Boolean algebra} with $1+1=1,$ and $B(n,0)=\mathbb{Z}_{n}$ for $%
n\geq 2.$ The semiring $B(n,i)$ is a semidomain if and only if $i\geq 1$ or $%
i=0$ and $n$ is a prime number. We call $B(n,i)$ with $i\geq 1$ the \textbf{%
Alarcon-Anderson semidomain}.
\end{ex}

\begin{punto}
Let $R$ be a (semi)ring and consider $Ideal(R),$ the complete lattice of all
ideals of $R$ with $I\vee J:=I+J$ and $I\wedge J:=I\cap J$ for ideals $I,J$
of $R.$ For $\emptyset \neq X\subseteq Ideal(R)\backslash \{R\},$ we say
that $R$ is an $X$\textbf{-top (semi)ring} iff $Ideal(R)$ is an $X$-top
lattice. For example, $R$ is a $Spec(R)$-top (semi)ring and the topology on $%
Spec(R)$ is the ordinary \emph{Zariski topology} on the spectrum of prime
ideals of $R$ (e.g., \cite{AM1969}). Moreover, for any $\emptyset \neq
Y\subseteq Spec(R),$ we have: $R$ is a $Y$-top (semi)ring by Corollary \ref%
{Y-top}. In particular, $R$ is a $Max(R)$-top (semi)ring and a $Min(R)$-top
(semi)ring, where $Max(R)$ (resp. $Min(R)$) is the spectrum of maximal
ideals (resp. minimal prime ideals) of $R.$
\end{punto}

\begin{defn}
(cf. \cite{Wis1991}, \cite{Lam2001}) We say that a commutative (semi)ring $R$
is \textbf{local} iff $R$ has a unique maximal ideal.
\end{defn}

\subsection*{General Topology}

In what follows, we recall some definitions and elementary results from
General Topology. We follow \cite{Wil1970} unless otherwise mentioned
explicitly.

\begin{notatation}
Let $X$ be a topological space. For $Y\subseteq X,$ we denote by $\mathcal{O}%
(Y)$ (resp. $\mathcal{C}(Y),$ $\mathcal{CO}(Y)$, $\mathcal{D}(Y),$ $\mathcal{%
K}(Y),$ $\mathcal{KO}(Y),$ $\mathcal{KC}(Y),$ $\mathcal{KCO}(Y),$ $\mathcal{I%
}(Y),$ $\mathcal{IC}(Y)$) the collection of all open (resp. closed, clopen,
connected, compact, compact open, compact closed, compact clopen,
irreducible closed, maximal irreducible closed) subsets of $X$ that contain $%
Y.$ By abuse of notation (compare with the notation for $Y=\emptyset $), we
denote by $\mathcal{O}(X)$ (resp. $\mathcal{C}(X),$ $\mathcal{CO}(X)$, $%
\mathcal{D}(X),$ $\mathcal{K}(X),$ $\mathcal{KO}(X),$ $\mathcal{KC}(X),$ $%
\mathcal{KCO}(X),$ $\mathcal{I}(X),$ $\mathcal{IC}(X)$) the collection of
all open (resp. closed, clopen, connected, compact, compact open, compact
closed, compact clopen, irreducible closed, maximal irreducible closed)
subsets of $X.$
\end{notatation}

\subsubsection*{Separation Axioms}

\begin{defns}
We say that a topological space $X$ is

\begin{enumerate}
\item \textbf{connected} iff $X\neq \emptyset $ and cannot be written as the
union of two \emph{disjoint proper} open (closed) subsets;

\item \textbf{hyperconnected} (\textbf{irreducible}) iff $X$ cannot be
written as the union of two \emph{proper} closed subsets, equivalently iff
no two non-empty \emph{open} subsets of $X$ are disjoint;

\item \textbf{ultraconnected} iff no two non-empty \emph{closed} subsets of $%
X$ are disjoint.
\end{enumerate}
\end{defns}

\begin{punto}
\label{anti-T2}Let $X$ be a topological space and set%
\begin{equation*}
\mathcal{S}(T_{2};X):=\{(x,y)\in (X\times X)\backslash \Delta (X)\mid
\exists \text{ }U\times V\in \mathcal{O}(x)\times \mathcal{O}(y)\text{ s.t. }%
U\cap V=\emptyset \}.
\end{equation*}%
We call $X$ \textbf{extremely non-Hausdorff} (resp. \textbf{anti-Hausdorff})
iff $\left\vert X\right\vert \geq 2$ and $\mathcal{S}(T_{2};X)=\emptyset $
(resp. $(Y\subseteq X$ is $T_{2}\Longrightarrow \left\vert Y\right\vert \leq
1)$).
\end{punto}

\begin{rem}
What we call \emph{extremely non-Hausdorff} spaces (as defined in \ref%
{anti-T2}) appeared in \cite[Theorem 4.2]{MM2012} under the name \emph{%
anti-Hausdorff} spaces. However, we reserve \textbf{anti-}$T_{2}$ as we
defined it to be consistent with the notion of \textbf{anti-}$\mathbf{P}$
spaces introduced by Bankston \cite{Ban1979}, which are roughly those
topological spaces that are \emph{almost hereditarily non}\textbf{-}$\mathbf{%
P}.$
\end{rem}

\begin{lem}
\label{anti-hyper}\emph{(cf. \cite[Theorem 4.2]{MM2012}, \cite{RV1980})} Let
$X$ be a topological space.

\begin{enumerate}
\item $X$ is extremely non-Hausdorff if and only if $\left\vert X\right\vert
\geq 2$ and $X$ is hyperconnected (irreducible).

\item $X$ is anti-Hausdorff if and only if $X$ is totally ordered.
\end{enumerate}
\end{lem}

\begin{punto}
Let $\mathcal{L}=(L;\vee ,0;\wedge ,1)$ be an $X$-top lattice for some $%
\emptyset \neq X\subseteq L\backslash \{1\}$ and consider the poset $(X,\leq
).$ We define the \textbf{height of} $x\in X$ as%
\begin{equation*}
ht(x):=\sup \{n\geq 0\text{ }\mid \text{there exists }\{x_{0},\cdots
,x_{n}\}\subseteq X\text{ with }x_{0}\lneqq \cdots \lneqq x_{n}=x\}.
\end{equation*}%
We define the \textbf{Krull dimension of }$X$ as $K.\dim (X):=\sup
\{ht(x)\mid x\in X\}.$
\end{punto}

\begin{punto}
Let $R$ be a commutative (semi)ring, $L=Ideal(R)$ and $X=Spec(R).$ The \emph{%
Krull dimension} $K.\dim (R)$ of $R$ is nothing but $K.\dim (Spec(R)).$ For
example, we have $K.\dim (\mathbb{Z})=K.\dim (Spec(\mathbb{Z}))=1$ and $%
K.\dim (\mathbb{W})=K.\dim (Spec(\mathbb{W}))=2.$
\end{punto}

\begin{defn}
(cf. \cite{G-R2007}, \cite{Dun1977}, \cite{Wil1970}) Let $X$ be a
topological space and set for $x\in X:$%
\begin{equation*}
Ker(x):=\bigcap\limits_{U\in \mathcal{O}(x)}U\text{ and }E(x):=\bigwedge%
\limits_{y\in X\backslash \{x\}}y.
\end{equation*}%
We say that $x\in X$ is \textbf{isolated} (resp. \textbf{kerneled}, \textbf{%
regular open}, \textbf{excluded}) iff $\{x\}$ is an open set (resp. $%
\{x\}=Ker(x),$ $\{x\}=int(\overline{\{x\}}),$ $E(x)=\bigwedge\limits_{y\in
D(x)}y$).
\end{defn}

For any topological space $X,$ we set%
\begin{equation*}
\begin{array}{lllllll}
K(X) & := & \{x\in X\text{ }\mid \text{ }\{x\}=Ker(\{x\})\}; &  & Iso(X) & :=
& \{x\in X\text{ }|\text{ }\{x\}\text{ is open}\}; \\
RO(X) & := & \{x\in X\text{ }|\text{ }\{x\}=int(\overline{\{x\}})\}; &  &
Cl(X) & := & \{x\in X\text{ }|\text{ }\{x\}\text{ is closed}\};%
\end{array}%
\end{equation*}

\begin{defn}
(e.g., \cite{PRV2009}, \cite{Dun1977}) A topological space $X$ is

\begin{enumerate}
\item $T_{\frac{1}{4}}$ iff any $x\in X$ is closed \emph{or} kerneled (i.e.,
iff $X=Cl(X)\cup K(X)$).

\item $T_{\frac{1}{2}}$ iff any $x\in X$ is closed \emph{or} isolated (i.e.,
iff $X=Cl(X)\cup Iso(X)$).

\item $T_{\frac{3}{4}}$ iff any $x\in X$ is closed \emph{or} regular open
(i.e., iff $X=Cl(X)\cup RO(X)$).
\end{enumerate}
\end{defn}

\begin{defn}
A topological space $X$ is said to be \textbf{sober} iff every irreducible
closed subset $Y\subseteq X$ has a \emph{unique} \textbf{generic point}
(i.e., $\exists !$ $y\in Y$ such that $Y=\overline{\{y\}}$).
\end{defn}

\begin{defn}
(cf. \cite{Est1988}, \cite{Hoc1969}, \cite{DST2019}) A \textbf{spectral space%
} is a topological space $X$ that satisfies any (hence all) of the following
equivalent conditions:

\begin{enumerate}
\item $X$ is sober, compact and has a base $\mathcal{B}\subseteq \mathcal{KO}%
(X)$ closed under finite intersections;

\item $X$ is homeomorphic to $Spec(R)$ for some \emph{commutative (semi)ring}
$R;$

\item $X$ is homeomorphic to a projective limit of \emph{finite }$T_{0}$
spaces.
\end{enumerate}
\end{defn}

\begin{defn}
A topological space is \textbf{quasi-Hausdorff} \cite{Hoc1969} iff for any $%
x\neq y$ in $X:$ \emph{either} $x$ and $y$ are separated by disjoint open
neighborhoods\emph{\ or }there exists $z\in X$ such that $\{x,y\}\subseteq
\overline{\{z\}}.$
\end{defn}

\begin{lem}
\label{T1-qH-T2}\emph{(\cite[Lemma 2.8.]{AA})}

\begin{enumerate}
\item Every finite $T_{0}$ space is spectral \emph{(cf. \cite{Hoc1969})}.

\item Every spectral space is quasi-Hausdorff \emph{(cf. \cite[Corollary 2,
page 45]{Hoc1969})}.

\item A topological space $X$ is $T_{2}$ if and only if $X$ is $T_{1}$ and
quasi-Hausdorff.

\item A $T_{1}$ spectral space is $T_{2}$ (cf.\emph{\ \cite[Exercise 3.11]%
{AM1969}}).
\end{enumerate}
\end{lem}

\begin{prop}
\label{dim-0}\emph{(\cite[Proposition 2.9.]{AA})} Let $\mathcal{L}=(L;\vee
,0;\wedge ,1)$ be an $X$-top lattice for some $\emptyset \neq X\subseteq
L\backslash \{1\}.$

\begin{enumerate}
\item $X$ is $T_{0}.$

\item If $X$ is finite, then $X$ is spectral.

\item $Max(X)=Max(C^{X}(\mathcal{L})).$

\item $X$ is $T_{1}\Longleftrightarrow $ $K.\dim (X)=0.$

\item $X$ is $T_{2}\Longleftrightarrow K.\dim (X)=0$ and $X$ is
quasi-Hausdorff.

\item $X$ is a $T_{\frac{1}{4}}\Longleftrightarrow K.\dim (X)\leq 1.$
\end{enumerate}
\end{prop}

\begin{thm}
(\cite{Lew1973})\ Let $X$ be a finite partially ordered set. Then there
exists a commutative ring $R$ such that $X\simeq Spec(R)$ (as posets).
\end{thm}

\section{Normal $X$-top lattices}

\qquad In this section, we focus on studying \emph{normal}, \emph{completely
normal} and \emph{perfectly normal} $X$-top lattices. We draw the attention
of the reader that when studying the aforementioned topological properties,
we do \textbf{not} assume the $T_{1}$ separation axiom as several references
do since $X$ is $T_{1}$ if and only if $K.\dim (X)=0$ (cf. Proposition \ref%
{dim-0} (4)). To avoid confusion and to make the manuscript self-contained,
we fix our terminology.

\begin{punto}
(cf. \cite{Wil1970}, \cite{RV1980}, \cite[Proposition 2]{GRV1981})\ Let $X$
be a topological space. Set%
\begin{equation*}
\begin{tabular}{lll}
$\mathcal{T}(\mathbf{R};X)$ & $:=$ & $\{(C,p)\text{ }\mid \text{ }\emptyset
\neq C\varsubsetneqq X\text{ is closed and }p\in X\backslash C\};$ \\
$\mathcal{S}(\mathbf{R};X)$ & $:=$ & $\{(C,p)\in \mathcal{T}(\mathbf{R}%
;X)\mid \exists $ $U\times V\in \mathcal{O}(C)\times \mathcal{O}(p)$ s.t. $%
U\cap V=\emptyset \};$ \\
$\mathcal{S}(\mathbf{CR};X)$ & $:=$ & $\{(C,p)\in \mathcal{T}(\mathbf{R}%
;X)\mid \exists \text{ }f:X\overset{\text{cts.}}{\longrightarrow }\mathbb{R}%
\text{ s.t. }f(C)=0\text{ and }f(p)=1\};$ \\
$\mathcal{T}(\mathbf{N};X)$ & $:=$ & $\{(C,D)$ $\mid $ $\emptyset \neq C,$ $%
D\subseteq X$ closed and $C\cap D=\emptyset \};$ \\
$\mathcal{S}(\mathbf{N};X)$ & $:=$ & $\{(C,D)\in \mathcal{T}(\mathbf{N}%
;X)\mid \exists $ $U\times V\in \mathcal{O}(C)\times \mathcal{O}(D)$ s.t. $%
U\cap V=\emptyset \}.$ \\
$\mathcal{T}(\mathbf{PN};X)$ & $:=$ & $\{(C,D)$ $\mid $ $C,$ $D\subseteq X$
closed and $C\cap D=\emptyset \};$ \\
$\mathcal{S}(\mathbf{PN};X)$ & $:=$ & $\{(C,D)\in \mathcal{T}(\mathbf{PN}%
;X)\mid \exists $ $f:X\overset{\text{cts.}}{\longrightarrow }\mathbb{R}\text{
s.t. }C=f^{-1}(0)\text{ and }D=f^{-1}(1)\}.$%
\end{tabular}%
\end{equation*}%
We say that $X$ is

\textbf{regular} (resp. \textbf{extremely non-regular, anti-regular}) iff $%
\mathcal{S}(\mathbf{R};X)=\mathcal{T}(\mathbf{R};X)$ (resp. $\mathcal{T}(%
\mathbf{R};X)\neq \emptyset $ and $\mathcal{S}(\mathbf{R};X)=\emptyset ,$ $%
Y\subseteq X$ is regular $\Longrightarrow \left\vert Y\right\vert \leq 1);$

\textbf{completely regular} iff $\mathcal{S}(\mathbf{CR};X)=\mathcal{T}(%
\mathbf{R};X);$

\textbf{normal} (resp. \textbf{extremely non-normal}, \textbf{anti-normal})\
iff $\mathcal{S}(\mathbf{N};X)=\mathcal{T}(\mathbf{N};X)$ (resp. $\mathcal{T}%
(\mathbf{N};X)\neq \emptyset $ and $\mathcal{S}(\mathbf{N};X)=\emptyset ,$ $%
Y\subseteq X$ is normal $\Longrightarrow \left\vert Y\right\vert \leq 2);$

\textbf{completely normal} (\textbf{hereditarily normal}) iff every subspace
$Y\subseteq X$ is normal;

\textbf{perfectly normal}\ iff $\mathcal{S}(\mathbf{PN};X)=\mathcal{T}(%
\mathbf{PN};X);$

$T_{3}$ (resp. $T_{3\frac{1}{2}},$ $T_{4},$ $T_{5},$ $T_{6}$) iff $X$ is $%
T_{1}$ and regular (resp. completely regular, normal, completely normal,
perfectly normal).
\end{punto}

\begin{remarks}
\label{e-non-regular-non-normal}Let $X$ be a topological space. It's clear
that:

\begin{enumerate}
\item If $X$ is extremely non-Hausdorff and $\mathcal{T}(\mathbf{R};X)\neq
\emptyset ,$ then $X$ is extremely non-regular.

\item If $X$ is extremely non-regular and $\mathcal{T}(\mathbf{N};X)\neq
\emptyset ,$ then $X$ is extremely non-normal.
\end{enumerate}
\end{remarks}

\begin{defn}
Let $X$ be a topological space. Two subsets $A,B\subseteq X$ are said to be
\textbf{separated} iff $A\cap \overline{B}=\emptyset =\overline{A}\cap B.$
\end{defn}

\begin{lem}
\label{CN-separated}(cf. \cite[15.B]{Wil1970}) A topological space $X$ is
completely normal if and only if for every pair $(A,B)\subseteq X\times X$
of \emph{separated} sets, there exists $U\times V\in \mathcal{O}(A)\times
\mathcal{O}(B)\ $such that $U\cap V=\emptyset .$
\end{lem}

\begin{punto}
Let $X$ be a topological space. A subset $A\subseteq X$ is said to be a $%
G_{\delta }$-set iff $A$ is a \emph{countable} intersection of open sets. We
say that $X$ is a $G_{\delta }$\textbf{-space} (a \textbf{perfect space}
\cite{Eng1989}) iff every closed set in $X$ is a $G_{\delta }$-set.
\end{punto}

\begin{lem}
\label{Vedenissoff}(e.g., \cite[Theorem 2, page 135]{Kur1966}) (The
Vedenissoff Theorem) A topological space $X$ is $X$ is perfectly normal if
and only if $X$ is a normal $G_{\delta }$-space.
\end{lem}

The proof of the following lemma is straightforward.

\begin{lem}
\label{hyp-e-non-regular}Let $X$ be a topological space.

\begin{enumerate}
\item If $X$ is ultraconnected, then $X$ is normal.

\item If $X$ is hyperconnected (irreducible) and $\mathcal{T}(\mathbf{R}%
;X)\neq \emptyset ,$ then $X$ is extremely non-regular.

\item If $X$ is extremely non-regular and $\mathcal{T}(\mathbf{N};X)\neq
\emptyset ,$ then $X$ is extremely non-normal.
\end{enumerate}
\end{lem}

\begin{remarks}
\label{rem-T3}(e.g., \cite{Kur1966}).

\begin{enumerate}
\item The Vedenissoff Theorem (cf. Lemma \ref{Vedenissoff}) is usually
stated for $T_{1}$ spaces (e.g. \cite[Theorem 1.5.19]{Eng1989}). However,
the proof does use this assumption.

\item For a topological space $X,$ we have the following implications which
are not reversible (in general):
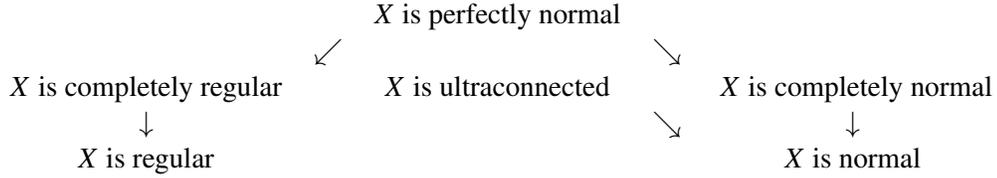
\begin{equation*}
\begin{tabular}{ccccc}
&  & $X\text{ is perfectly normal}$ &  &  \\
& $\swarrow $ &  & $\searrow $ &  \\
$X\text{ is completely }$regular &  & $X\text{ is ultraconnected}$ &  & $%
\text{ }X\text{ is completely }$normal \\
$\downarrow $ &  &  & $\searrow $ & $\downarrow $ \\
$X\text{ is regular}$ &  &  &  & $X\text{ is }$normal%
\end{tabular}%
\end{equation*}

\item Every $T_{0}$ regular (resp. completely regular, perfectly normal)
space is $T_{1}$, whence $T_{3}$ (resp. $T_{3\frac{1}{2}},$ $T_{6}$).

\item We have%
\begin{equation*}
T_{6}\Longrightarrow T_{5}\Longrightarrow T_{4}\Longrightarrow T_{3\frac{1}{2%
}}\Longrightarrow T_{3}\Longrightarrow T_{2}\Longrightarrow
T_{1}\Longrightarrow T_{\frac{3}{4}}\Longrightarrow T_{\frac{1}{2}%
}\Longrightarrow T_{\frac{1}{4}}\Longrightarrow T_{0}.
\end{equation*}
\end{enumerate}
\end{remarks}

\begin{punto}
Let $\mathcal{L}=(L;\vee ,0;\wedge ,1)$ be an $X$-top lattice for some $%
\emptyset \neq X\subseteq L\backslash \{1\}.$ We say that

$X$ has the \textbf{pm-property} iff $\left\vert Max(x;X)\right\vert =1$ for
every $x\in X;$

$X$ has the \textbf{m-property} iff $\left\vert Min(x;X)\right\vert =1$ for
every $x\in X;$

$X$ is\textbf{\ Jacobson} iff for every $x\in X,$ we have $%
x=\bigwedge\limits_{Max(x;X)}\mathfrak{m};$

$X$ is\textbf{\ dual Jacobson }iff for every $x\in X,$ we have $%
x=\bigvee\limits_{Min(x;X)}m.$

We say that a (semi)ring $R$ is a \textbf{pm-(semi)ring} \cite{MO1971}
(resp. an \textbf{m}-(\textbf{semi})\textbf{ring} \cite{Avi2006}, a \textbf{%
Jacobson }(\textbf{semi})\textbf{ring,} a \textbf{dual Jacobson }(\textbf{%
semi})\textbf{ring)} iff $Spec(R)\ $has the pm-property (resp. has the
m-property, is Jacobson, is dual Jacobson). For more on Jacobson commutative
rings consult \cite{AM1969}.
\end{punto}

\begin{rem}
The pm-(semi)rings were called \emph{Gelfand }(\emph{semi})\emph{rings} in
\cite[Theorem 4.2]{PRV2009}. We choose not to use this terminology to avoid
any possible confusion with \emph{Gelfand semirings }in the sense of \cite%
{Gol1999}.
\end{rem}

De Marco and Orsatti in \cite{MO1971} characterized the commutative rings
for which the prime spectrum is \emph{normal} (not necessarily $T_{1}$)\ as
the pm-rings. For $X$-top lattices, we investigate the relation between the
normality of $X$ and the pm-property.

\begin{defn}
Let $\mathcal{L}=(L;\vee ,0;\wedge ,1)$ be an $X$-top lattice for some $%
\emptyset \neq X\subseteq L\backslash \{1\}.$ A subspace $Y\subseteq X$ is
said to be a \textbf{retract} of $X$ iff there exists a continuous map
(called \textbf{a retraction}) $f:X\longrightarrow Y$ such that $f_{\mid
_{Y}}=id_{Y}.$
\end{defn}

The following is a restatement of \cite[Theorem 2.1]{MO1971} and \cite[%
Theorem 4.2]{PRV2009} (cf. \cite{Hoc1969}):

\begin{prop}
\label{tspm}Let $\mathcal{L}=(L;\vee ,0;\wedge ,1)$ be an $X$-top lattice
for some $\emptyset \neq X\subseteq L\backslash \{1\}.$ If $X$ is spectral,
then the following are equivalent:

\begin{enumerate}
\item $X$ has the pm-property;

\item $Max(X)$ is a retract of $X$;

\item $X$ is a normal space;

\item $X$ is homeomorphic to $Spec(R)$ for some $pm$-(semi)ring.
\end{enumerate}
\end{prop}

\begin{rem}
\label{ret-unique}Let $\mathcal{L}=(L;\vee ,0;\wedge ,1)$ be an $X$-top
lattice for some $\emptyset \neq X\subseteq L\backslash \{1\}.$ Assume that $%
X$ is max-retractable with retraction $\mathfrak{r}:X\longrightarrow Max(X).$
Let $x\in X$ and set $\mathfrak{m=r}(x).$ Notice that $Max(X)$ is $T_{1},$
whence $\{\mathfrak{m}\}$ is closed in $Max(X).$ Since $\mathfrak{r}$ is
continuous, $\mathfrak{r}^{-1}(\mathfrak{m})$ is closed in $Spec(X).$ It
follows that $V(x):=\overline{\{x\}}\subseteq \mathfrak{r}^{-1}(\mathfrak{m}%
),$ i.e., $\mathfrak{r}(V(x))\subseteq \{\mathfrak{m}\}.\blacksquare $
\end{rem}

The following result investigates to which extent, Proposition \ref{tspm}
can be generalized to $X$-top lattices for which $X$ is \emph{not spectral}.

\begin{prop}
\label{txnpm}Let $\mathcal{L}=(L;\vee ,0;\wedge ,1)$ be an $X$-top lattice
for some $\emptyset \neq X\subseteq L\backslash \{1\}$ and assume that $X$
is coatomic.

\begin{enumerate}
\item If $X$ is normal, then $X$ has the pm-property.

\item If $Max(X)$ is a retract of $X,$ then $X$ has the pm-property.
\end{enumerate}
\end{prop}

\begin{Beweis}
Let $X$ be coatomic.

\begin{enumerate}
\item Assume that $X$ is normal. Let $x\in X$ and notice that $Max(x;X)\neq
\emptyset $ (since $X$ is coatomic). Suppose that there exists $x\in X$ and $%
\mathfrak{m}\neq \mathfrak{m}^{\prime }$ in $Max(x;X).$ Notice that $\{%
\mathfrak{m}\}=V(\mathfrak{m})$ and $\{\mathfrak{m}^{\prime }\}=V(\mathfrak{m%
}^{\prime })$ are closed sets in $X$ but cannot be separated by disjoint
open sets: any open set $D(a)$ (where $a\in L$) that contains $\{\mathfrak{m}%
\}$ or $\{\mathfrak{m}^{\prime }\}$ will contain $x$ as well, a
contradiction.

\item Assume that $\mathfrak{r}:X\longrightarrow Max(X)$ is a \emph{%
retraction.} Let $x\in X$ and consider $\mathfrak{m}:=\mathfrak{r}(x).$ Then
$\mathfrak{r}(Max(x;X))\subseteq \mathfrak{r}(V(x))\overset{\text{Remark \ref%
{ret-unique}}}{=}\{\mathfrak{m}\}.$ Since $X$ is coatomic we have $%
Max(x;X)\neq \emptyset .$ Since $f_{\mid _{Max(X)}}=id_{Max(X)},$ it follows
by Remark \ref{ret-unique} that $Max(x;X)=\{\mathfrak{m}\}.$ Since $x\in X$
was arbitrary, we conclude that $X$ has the pm-property.$\blacksquare $
\end{enumerate}
\end{Beweis}

The following proposition gives two characterizations of \emph{normal} $X$%
-top lattices under some conditions.

\begin{thm}
\label{pxn}Let $\mathcal{L}=(L;\vee ,0;\wedge ,1)$ be an $X$-top lattice for
some $\emptyset \neq X\subseteq L\backslash \{1\}$ and assume that $X$ is
coatomic.

\begin{enumerate}
\item Assume that $X$ is atomic and $Min(X)$ is finite. Then $X$ is normal
if and only if $X$ has the pm-property.

\item Assume that $X$ is \emph{completely} strongly $X$-irreducible and $%
Max(X)$ is finite. Then $X$ is a retract of $Max(X)$ if and only if $X$ has
the pm-property.
\end{enumerate}
\end{thm}

\begin{Beweis}
\begin{enumerate}
\item If $X$ is coatomic and normal, then $X$ has the pm-property by Theorem %
\ref{txnpm} (1).

For the converse, assume that $X$ has the pm-property, $X$ is atomic and $%
|Min(X)|<\infty .$ Then $X$ is clearly coatomic. Since $X$ has the
pm-property, there is a well-defined map%
\begin{equation*}
\mu :X\longrightarrow Max(X),\text{ where }x\leq \mu (x).
\end{equation*}%
Let $\emptyset \neq V(c),\ V(d)\varsubsetneqq X$ be \emph{disjoint.} Since $%
X $ is atomic, we can pick for each $x\in V(c)\cup V(x)$ some $a_{x}\in
Min(X)$ such that $a_{x}\leq x\leq \mu (x)$ (notice that $\mu (a_{x})=\mu
(x) $ since $X$ has the pm-property). Set%
\begin{equation*}
\begin{array}{ccccccc}
J & := & Max(c;X); &  & K & := & Max(d;X); \\
\widetilde{J} & := & \{a\in Min(X)\mid \mu (a)\notin J\}; &  & \widetilde{K}
& := & \{a\in Min(X)\mid \mu (a)\notin K\}; \\
U & := & \bigcap\limits_{a\in \widetilde{J}}D(a), &  & V & := &
\bigcap\limits_{a\in \widetilde{K}}D(a).%
\end{array}%
\end{equation*}%
Notice that $U$ and $K$ are open in $X$ (since $Min(X)$ is finite).
Moreover, $V(c)\subseteq U:$ Let $x\in V(c)$ and suppose that $x\notin U.$
Then there exists $a_{x}\in Min(X)$ such that $a_{x}\leq x$ with $\mu
(a_{x})\in Max(X)\backslash J.$ Since $X$ has the pm-property, it follows
that $\mu (x)=\mu (a_{x})\notin V(c),$ a contradiction. One can prove,
similarly, that $V(d)\subseteq V.$ Notice that $J\cap K=\emptyset ,$ whence $%
Min(X)=\widetilde{J}\cup \widetilde{K}$ and it follows that $U\cap
V=\emptyset .$ So, $X$ is normal.

\item If $Max(X)$ is a retract of $X,$ then $X$ has the pm-property by
Theorem \ref{txnpm} (2).

For the converse, assume that $X$ has the pm-property, $X$ is \emph{%
completely} strongly $X$-irreducible and $Max(X)$ is finite. The map $%
\mathfrak{r}{:X\longrightarrow Max(X)}$ defined by $\mathfrak{r}(x)=%
\mathfrak{m}$, where $\mathfrak{m}$ is the \emph{unique} element in $%
Max(x;X):=V(x)\cap Max(X)$ is well defined and satisfies $\mathfrak{r}_{\mid
_{Max(X)}}=id_{Max(X)}.$

Consider $U_{X}(\mathfrak{m}):=\{z\in X\ |\ z\leq \mathfrak{m}\}.$ \textbf{%
Claim:} $V(\bigwedge\limits_{z\in U_{X}(\mathfrak{m})}z)=U_{X}(\mathfrak{m}%
). $

Suppose that there exists $y\in V(\bigwedge\limits_{z\in U_{X}(\mathfrak{m}%
)}z)\backslash U_{X}(\mathfrak{m}).$ Since $X$ is coatomic, there exists $%
\mathfrak{m}^{\prime }\in Max(X)\backslash \{\mathfrak{m}\}$ such that $%
y\leq \mathfrak{m}^{\prime }$. By assumption, $X$ is \emph{completely}
strongly $X$-irreducible and so there exists $z_{y}\in U_{X}(\mathfrak{m})$
such that $z_{y}\leq y.$ It follows that $\{\mathfrak{m},\mathfrak{m}%
^{\prime }\}\subseteq Max(z_{y};X)$, a contradiction (to the assumption that
$X$ has the pm-property). So, $V(\bigwedge\limits_{z\in U_{X}(\mathfrak{m}%
)}z)\subseteq U_{X}(\mathfrak{m}).$ The reverse inclusion is trivial.

So, we have%
\begin{equation*}
\mathfrak{r}^{-1}(V(\mathfrak{m}))=\mathfrak{r}^{-1}(\mathfrak{m})=\{z\in X\
|\ z\leq \mathfrak{m}\}=U_{X}(\mathfrak{m})=V(\bigwedge\limits_{z\in U_{X}(%
\mathfrak{m})}z).
\end{equation*}%
Since $Max(X)$ is finite, we conclude that $\mathfrak{r}$ is a continuous
map. So, $Max(X)$ is a retract of $X.\blacksquare $
\end{enumerate}
\end{Beweis}

\begin{thm}
\label{FMax-pm-retract}Let $\mathcal{L}=(L;\vee ,0;\wedge ,1)$ be an $X$-top
lattice for some $\emptyset \neq X\subseteq L\backslash \{1\}.$ If $X$ is
coatomic, atomic with both $Min(X)$ and $Max(X)$ finite, then%
\begin{equation*}
X\text{ is normal }\Longleftrightarrow \text{ }X\text{ has the pm-property }%
\Longleftrightarrow \text{ }Max(X)\text{ is a retract of }X.
\end{equation*}
\end{thm}

The following result is a direct consequence of Theorem \ref{FMax-pm-retract}
as well as a direct consequence of Proposition \ref{tspm} since any finite $%
T_{0}$ spaces is spectral (cf. Lemma \ref{T1-qH-T2}).

\begin{cor}
\label{finite-pm-retract}Let $L$ is an $X$-top lattice for some $\emptyset
\neq X\subseteq L\backslash \{1\}.$ If $X$ is \emph{finite}, then the
following are equivalent:

\begin{enumerate}
\item $Max(X)$ is a retract of $X;$

\item $X$ has the pm-property;

\item $X$ is a normal space;

\item $X$ is homeomorphic to $Spec(R)$ for some pm-(semi)ring with finitely
many primes.
\end{enumerate}
\end{cor}

\begin{prop}
\label{normal-local}Let $\mathcal{L}=(L;\vee ,0;\wedge ,1)$ be an $X$-top
lattice for some $\emptyset \neq X\subseteq L\backslash \{1\}.$

\begin{enumerate}
\item If $X$ is local, then $X$ is normal.

\item If $X$ is colocal and normal (and coatomic), then $\left\vert
Max(X)\right\vert \leq 1$ ($X$ is local).
\end{enumerate}
\end{prop}

\begin{Beweis}
\begin{enumerate}
\item Let $X$ be local (i.e., $X$ is coatomic and $Max(X)=\{\mathfrak{m}\}$%
). Notice that $\mathfrak{m}\in V(a)$ for every $a\in L.$ Consequently, $X$
is ultraconnected and so \emph{trivially} normal.

\item Assume that $X$ is normal, colocal with $Min(X)=\{m\}$ and $Max(X)\neq
\emptyset .$ Suppose that there exist $\mathfrak{m}\neq \mathfrak{m}^{\prime
}$ in $Max(X).$ Notice that $\{\mathfrak{m},\mathfrak{m}^{\prime
}\}\subseteq V(m).$ The singletons $\{\mathfrak{m}\}=V(\mathfrak{m})$ and $\{%
\mathfrak{m}^{\prime }\}=V(\mathfrak{m}^{\prime })$ are closed but cannot be
separated by disjoint open disjoint sets: any open set $D(a)$ (where $a\in L$%
) that contains $\{\mathfrak{m}\}$ or $\{\mathfrak{m}^{\prime }\}$ would
contain $m$ as well, \emph{i.e.,} one cannot find \emph{disjoint} open sets
separating $\{\mathfrak{m}\}$ and $\{\mathfrak{m}^{\prime }\},$ a
contradiction. So, $\left\vert Max(X)\right\vert =1.\blacksquare $
\end{enumerate}
\end{Beweis}

\begin{thm}
\label{cor-local}Let $\mathcal{L}=(L;\vee ,0;\wedge ,1)$ be an $X$-top
lattice for some $\emptyset \neq X\subseteq L\backslash \{1\}.$ If $X$ is
colocal and coatomic, then the following are equivalent:

\begin{enumerate}
\item $X$ is local;

\item $X$ is ultraconnected;

\item $X$ is normal.
\end{enumerate}
\end{thm}

\begin{ex}
\label{local-semi(ring)}Let $R$ be a local (semi)ring. Then $Spec(R)\ $is
ultraconnected, whence normal.
\end{ex}

\begin{cor}
\label{domain-normal}Let $R$ be a (semi)domain. The following are equivalent:

\begin{enumerate}
\item $R$ is local;

\item $Spec(R)$ is ultraconnected;

\item $Spec(R)$ is normal;

\item $Spec(R)$ is max-retractable;

\item $R$ is a pm-(semi)ring.
\end{enumerate}
\end{cor}

\begin{ex}
\label{pm-not-normal}The integral domain $\mathbb{Z}$ is \emph{not} local,
whence $X=Spec(\mathbb{Z})$ is \emph{not} normal by Corollary \ref%
{domain-normal}. In fact, $X$ is hyperconnected (irreducible): for any
positive integers $m\neq n,$ we have%
\begin{equation*}
D(m\mathbb{Z})\cap D(n\mathbb{Z})=D(mn\mathbb{Z})=\{p\mathbb{Z}\mid p\nmid
mn\}\neq \emptyset .
\end{equation*}%
Moreover, $\left\vert \mathcal{T}(\mathbf{R};Spec(\mathbb{Z}))\right\vert
=\infty ,$ since $p\mathbb{Z}\in X\backslash V(m\mathbb{Z})$ (where $p\in
\mathbb{P}$) if and only if $p\nmid m.$ Moreover, $\left\vert \mathcal{T}(%
\mathbf{N};Spec(\mathbb{Z}))\right\vert =\infty ,$ since $V(m\mathbb{Z})\cap
V(n\mathbb{Z})=\emptyset $ if and only if $g.c.d.(m,n)=1.$ So, $X$ is
extremely non-Hausdorff, whence extremely non-regular and consequently
extremely non-normal by Lemma \ref{e-non-regular-non-normal}.$\blacksquare $
\end{ex}

\begin{ex}
Consider the semiring $S=(\mathbb{W}\cup \{\infty \};\oplus ,0;\oplus
,\infty ),$ where%
\begin{equation*}
a\oplus b=\max \{a,b\}\text{ and }a\oplus b=\min \{a,b\}.
\end{equation*}%
Setting $J_{s}=\{x\in S\mid x\leq s\},$ we have $J_{0}=\{0\}$ and $J_{\infty
}=S.$ Notice that $Spec(S)=\{J_{s}$ $\mid $ $s\in \mathbb{W}\}\cup \{\mathbb{%
W}\}.$

Notice that $S$ is local with $Max(S)=\{\mathbb{W}\}$ and colocal with $%
Min(S)=\{0\}.$ It follows, by Theorem \ref{cor-local}, that $Spec(S)$ is
(perfectly, completely) normal. Clearly, $S$ is a pm-semiring and%
\begin{equation*}
\mu :Spec(S)\longrightarrow Max(S),\text{ }J_{i}\longmapsto \mathbb{W}
\end{equation*}%
is a retraction.

Let $Y:=Spec(S)\backslash \{\mathbb{W}\}:=\{J_{s}\mid s\in \mathbb{W}\}.$
Then $Y$ is clearly ultraconnected, whence normal. However, $Y$ is neither
max-retractable, nor has the pm-property, nor local (as $Max(Y)=\emptyset $%
). Notice that $Y$ is atomic, colocal with $Min(Y)=\{0\}$ and both of $%
Min(Y) $ and $Max(Y)=\emptyset $ are finite. This shows that the assumption
that $Y$ is \emph{coatomic} cannot be dropped from the assumptions of
Proposition \ref{txnpm} (1), Theorem \ref{pxn} (1) and Theorem \ref%
{FMax-pm-retract}. Moreover, it shows that the finiteness condition cannot
be dropped from the assumption of Corollary \ref{finite-pm-retract}.$%
\blacksquare $
\end{ex}

We recall a useful description of the prime spectra of the semidomains $%
B(n,i):$

\begin{punto}
\label{Spec(B(n,i))}(\cite[Theorem 24]{AA1994}) Let $n\geq 2,$ $1\leq i\leq
n-1$ and set%
\begin{equation*}
\mathfrak{m}_{n}:=\{0,2,3,...,n-1\}\text{ for }n\geq 3).
\end{equation*}

\begin{enumerate}
\item $K.\dim (B(n,i))=0$ if $i=0$ or $n=2$ and $i=1.$

\item $K.\dim (B(n,i))=1$ if $n\geq 3$ and $i=1$. In this case,%
\begin{equation*}
Spec(B(n,i)=\{0\}\cup \{pB(n,i)\text{ }|\text{ }p\text{ is a prime divisor
of }n-1\}.
\end{equation*}

\item $K.\dim (B(n,i))=1$ if $n\geq 3$ and $i=n-1$. In this case, $%
Spec(B(n,i)=\{0,\mathfrak{m}_{n}\}.$

\item $K.\dim (B(n,i))=2$ if $n\geq 4$ and $2\leq i\leq n-2$. In this case,
\begin{equation*}
Spec(B(n,i)=\{0,\mathfrak{m}_{n}\}\cup \{pB(n,i)\text{ }|\text{ }p\text{ is
a prime divisor of }n-i\}.
\end{equation*}
\end{enumerate}
\end{punto}

\begin{ex}
\label{local-ultra}Consider the pentagon $N_{5}:$
\begin{figure}[]
\centering
\begin{equation*}
\xymatrix{ & 1 {\bullet} \ar@{-}[dl] \ar@{-}[ddr] & \\ x { \bullet}
\ar@{-}[d] & & \\ y { \bullet} \ar@{-}[dr] & & { \bullet} z \ar@{-}[dl] \\ &
0 {\bullet} &}
\end{equation*}%
\caption{$N_5$: A non-distributive $X$-top lattice}
\label{N_5}
\end{figure}
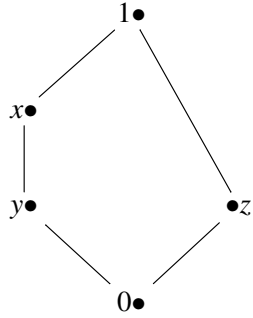
Although $N_{5}$ is \emph{not} distributive (cf. \cite[Theorem 101]{Gra2010}%
), we can still find $X$ such that $N_{5}$ is an $X$-top lattice.

\begin{enumerate}
\item Let $X=\{0,y,x\}=C^{X}(N_{5}).$ Clearly $X=SI^{C^{X}(N_{5})}(X),$
whence $N_{5}$ is an $X$-top lattice by Theorem \ref{xct}. The collection of
closed sets is%
\begin{equation*}
\begin{tabular}{lllllllllll}
$V_{X}(0)$ & $=$ & $X,$ &  & $V_{X}(y)$ & $=$ & $\{x,y\},$ &  & $V_{X}(z)$ &
$=$ & $\emptyset ,$ \\
$V_{X}(1)$ & $=$ & $\emptyset ,$ &  & $V_{X}(x)$ & $=$ & $\{x\}.$ &  &  &  &
\end{tabular}%
\end{equation*}%
Notice that $X$ is colocal with $Min(X)=\{0\}$ and local with $Max(X)=\{x\}.$
Clearly, $X$ is ultraconnected, whence $X$ is \emph{trivially} normal.

On the other hand, the collection of open sets of $X$ is given by%
\begin{equation*}
\begin{tabular}{lllllllllll}
$D_{X}(0)$ & $=$ & $\emptyset ,$ &  & $D_{X}(y)$ & $=$ & $\{0\},$ &  & $%
D_{X}(z)$ & $=$ & $X,$ \\
$D_{X}(1)$ & $=$ & $X,$ &  & $D_{X}(x)$ & $=$ & $\{0,y\}.$ &  &  &  &
\end{tabular}%
\end{equation*}%
Clearly, $X$ is hyperconnected (irreducible), whence $X$ is extremely
non-Hausdorff by Lemma \ref{anti-hyper}. Since $\mathcal{T}(\mathbf{R}%
,X)\neq \emptyset ,$ we conclude that $X$ is extremely non-regular by Lemma %
\ref{hyp-e-non-regular}. Notice that $X$ is a chain $0\lvertneqq y\lvertneqq
x$ and $K.\dim (X)=2,$ whence $X$ is not $T_{1}$ (cf. Proposition \ref{dim-0}
(4)) In fact, $X$ is not even $T_{\frac{1}{4}}$ by Proposition \ref{dim-0}
(6).

Since $X$ is finite and $T_{0},$ it turns out that $X$ is spectral and $%
X\approx Spec(R),$ where $R$ is any \emph{valuation ring} with $K.\dim (R)=2$
(e.g., $R=\mathbb{R}[[x,y]]+y\cdot \mathbb{R}((x))[[y]]$ since $%
Spec(R)=\{(0),$ $(y),(x,y)$). Moreover, $X$ is also homeomorphic to $%
Spec(B(n,i)),$ where $n\geq 4$ and $2\leq i\leq n-2$ with $n-i$ prime (cf.
Proposition \ref{Spec(B(n,i))} (4)). In particular, $X\approx
Spec(B(6,3))=\{0,\{0,3\},\{0,2,3,4,5\}\}.$

\item Let $Y=\{y,x\}\subseteq X.$ Then $C^{Y}(N_{5})=\{y,x\}=Y.$ By Then $%
N_{5}$ is a $Y$-top lattice by and (1). The collection of closed sets is%
\begin{equation*}
\begin{tabular}{lllllllllll}
$V_{Y}(0)$ & $=$ & $Y,$ &  & $V_{Y}(y)$ & $=$ & $Y,$ &  & $V_{Y}(z)$ & $=$ &
$\emptyset ,$ \\
$V_{Y}(1)$ & $=$ & $\emptyset ,$ &  & $V_{Y}(x)$ & $=$ & $\{x\}.$ &  &  &  &
\end{tabular}%
\end{equation*}%
Notice that $Y$ is colocal with $Min(Y)=\{y\}$ and local with $Max(Y)=\{x\}.$
Clearly, $Y$ is ultraconnected, whence $X$ is \emph{trivially} normal.

Since $Y$ is finite and $T_{0},$ we know $Y$ is spectral. In fact, $Y\approx
Spec(R),$ where $R$ is any DVR (e.g. $R$ is the discrete valuation ring $%
\mathbb{W}[[x]]$). Moreover, $X\approx Spec(B(n,n-1))=\{\{0\},\{0,2,\cdots
,n-1\}\}$ for any $n\geq 3$ (cf. Proposition \ref{Spec(B(n,i))} (3)).

On the other hand, the collection of open sets of $Y$ is given by%
\begin{equation*}
\begin{tabular}{lllllllllll}
$D_{Y}(0)$ & $=$ & $\emptyset ,$ &  & $D_{Y}(y)$ & $=$ & $\emptyset ,$ &  & $%
D_{Y}(z)$ & $=$ & $Y,$ \\
$D_{Y}(1)$ & $=$ & $Y,$ &  & $D_{Y}(x)$ & $=$ & $\{y\}.$ &  &  &  &
\end{tabular}%
\end{equation*}%
Clearly, $X$ is hyperconnected (irreducible), whence $X$ is extremely
non-Hausdorff by Lemma \ref{anti-hyper}. Since $\mathcal{T}(\mathbf{R}%
,X)\neq \emptyset ,$ we conclude that $X$ is extremely non-regular by Lemma %
\ref{hyp-e-non-regular}. Notice that $x$ is closed while $y$ is isolated but
not regular open, whence $Y$ is $T_{\frac{1}{2}}$ but not $T_{3\frac{1}{4}}$%
. $\blacksquare $
\end{enumerate}
\end{ex}

\begin{ex}
The semiring $\mathbb{W}$ is local with maximal ideal $\mathbb{W}\backslash
\{1\}$ and its prime spectrum is
\begin{figure}[]
\centering
\begin{equation*}
\xymatrix{& & & {\mathbb{W}}\backslash{{\{1\}}} \ar@{-}[ddll] \ar@{-}[ddl]
\ar@{-}[dd] \ar@{-}[ddr] \ar@{-}[ddr] \ar@{-}[ddlll] & \\ & & & & \\
2{\mathbb{W}} \ar@{-}[ddrrr] & 3{\mathbb{W}} \ar@{-}[ddrr] & \cdots
\ar@{-}[ddr] & p {\mathbb{W}} \ar@{-}[dd] & \cdots \ar@{-}[ddl] \\ & & & &
\\ & & & 0 & }
\end{equation*}%
\caption{The prime spectrum of ${\mathbb{W}}$}
\label{Spec(W)}
\end{figure}
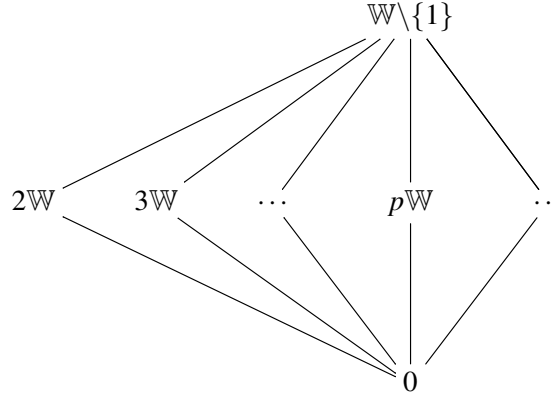
Notice that $Min(\mathbb{W})=0$ and $Max(\mathbb{W})=\{\mathbb{W}\backslash
\{1\}\},$ whence $X$ is ultraconnected and consequently \emph{trivially}
normal. Notice that $X$ is \emph{not} perfectly normal as one might think.
For example, the closed set $V(\mathbb{W}\backslash \{1\})=\{\mathbb{W}%
\backslash \{1\}\}$ is \emph{not} a $G_{\delta }$-set: if $\mathbb{W}%
\backslash \{1\}\in D(I)$ for some ideal $I,$ then $I\nsubseteqq \mathbb{W}%
\backslash \{1\},$ a contradiction.

Consider $Y:=X\backslash \{\mathbb{W}\backslash \{1\}\}\subseteq Spec(%
\mathbb{W}).$ Then $\mathbb{W}$ is a $Y$-top semiring by Corollary \ref%
{Y-top}. Moreover, $Y=D(\mathbb{W}\backslash \{1\})$ is an \emph{open}
subspace of $X.$ Clearly, $Y\approx Spec(\mathbb{Z}).$ So, $Y$ is
hyperconnected (irreducible), whence extremely non-Hausdorff, extremely
non-regular and extremely non-normal!$\blacksquare $
\end{ex}

Before we proceed, we consider a special class of graphs that will be used
throughout the rest of the section. In Graph Theory, a \textbf{tree} if the
a connected acyclic undirected graph, i.e., an undirected graph in which
every pair of distinct vertices is connected by \emph{exactly} one path.
However, we are interested in a very special class of \emph{rooted trees}
associated to posets.

\begin{punto}
Let $(P,\leq )$ be a partially ordered set.

\begin{enumerate}
\item If $(P,\leq )$ is a chain of $n+1$ elements (abbreviated $\mathcal{C}%
_{n+1}$), then we say that $P$ is a \textbf{chain of length }$n$.

\item A $\bigwedge $-\textbf{tree} is a non-empty subset $\mathcal{T}%
\subseteq P$ that satisfies the following conditions $\forall x,y,z\in
\mathcal{T}$:

\begin{enumerate}
\item if $x||y\ \in \mathcal{T}$, i.e., $x$ and $y$ are \emph{incomparable},
then $\exists $ $z\in \mathcal{T}$ such that $z\gvertneqq x$ and $%
z\gvertneqq y$;

\item if $x\lvertneqq y$ and $x\lvertneqq z$, then $y$ and $z$ are \emph{%
comparable}.
\end{enumerate}

If $\mathcal{T}$ is a $\bigwedge $-tree and $Min(\mathcal{T})$ is finite,
then we say that $\mathcal{T}$ is of a \textbf{finite base} (in this case, $%
\mathcal{T}$ is necessarily finite and has a \emph{unique} maximal element $%
\mathfrak{m}$ and we denote it by $(\mathcal{T},\mathfrak{m})$). We denote
by $\mathcal{T}_{n}$ the $\bigwedge $-tree of \emph{height} $1$ and $%
\left\vert Min(\mathcal{T})\right\vert =n.$ A collection of \emph{disjoint} $%
\bigwedge $-trees is called a $\bigwedge $-\textbf{forest}.

With a $\bigvee $\textbf{-tree} $\mathcal{V}$, we mean a $\bigwedge $-tree
in the dual poset $P^{o}=(P,\geq ).$ If $\mathcal{V}$ is a $\bigvee $-tree
with $Max(\mathcal{V})$ finite, then we say that $\mathcal{V}$ is of a
\textbf{finite cover} (in this case, $\mathcal{V}$ is necessarily finite and
has a \emph{unique} minimal element $m$ and we denote it by $(\mathcal{V},m)$%
). We denote by $\mathcal{V}_{n}$ the $\bigvee $-tree of \emph{height} $1$
and $\left\vert Max(\mathcal{V})\right\vert =n.$

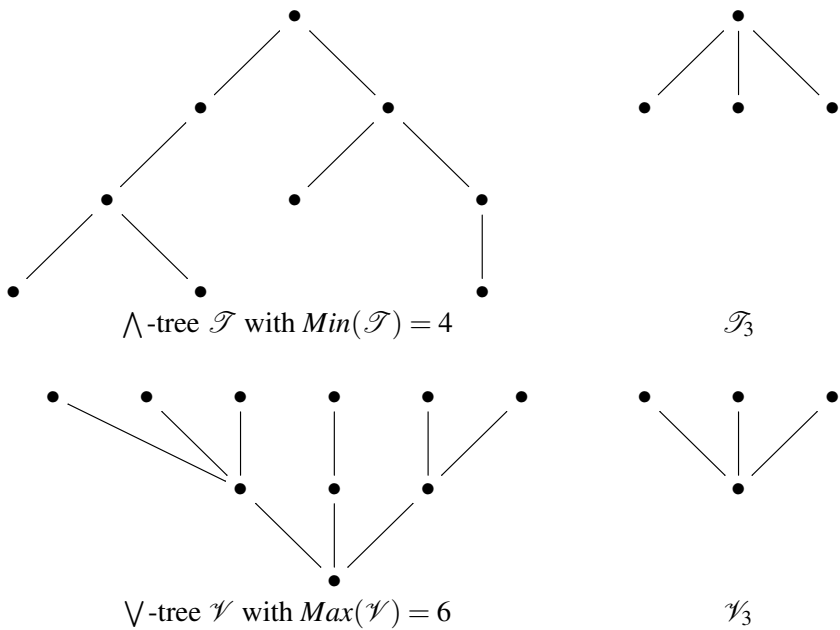
\begin{figure}[]
\centering
\begin{equation*}
\begin{array}{ccc}
\xymatrix{& & & \bullet \ar@{-}[dr] \ar@{-}[dl] & & & \\ & & \bullet
\ar@{-}[dl] & & \bullet \ar@{-}[dl] \ar@{-}[dr] & & \\ & \bullet \ar@{-}[dl]
\ar@{-}[dr] & & \bullet & & \bullet \ar@{-}[d] & \\ \bullet & & \bullet & &
& \bullet & } &  & \xymatrix{& \bullet \ar@{-}[dl] \ar@{-}[d] \ar@{-}[dr] &
\\ \bullet & \bullet & \bullet} \\
\bigwedge \text{-tree }\mathcal{T}\text{ with }Min(\mathcal{T})=4 &  &
\mathcal{T}_{3} \\
&  &  \\
\xymatrix{\bullet \ar@{-}[drr] & \bullet \ar@{-}[dr] & \bullet \ar@{-}[d] &
\bullet \ar@{-}[d] & \bullet \ar@{-}[d] & \bullet \ar@{-}[dl] \\ & & \bullet
& \bullet & \bullet & \\ & & & \bullet \ar@{-}[ul] \ar@{-}[u] \ar@{-}[ur] &
& } &  & \xymatrix{\bullet & \bullet & \bullet \\ & \bullet \ar@{-}[ul]
\ar@{-}[u] \ar@{-}[ur] } \\
\bigvee \text{-tree }\mathcal{V}\text{ with }Max(\mathcal{V})=6 &  &
\mathcal{V}_{3}%
\end{array}%
\end{equation*}%
\caption{Examples of trees}
\label{trees}
\end{figure}

Notice that for the chain $\mathcal{C}_{2},$ we have $\mathcal{T}_{1}=%
\mathcal{C}_{2}=\mathcal{V}_{1},$ i.e. $\mathcal{C}_{2}$ is a $\bigwedge $%
-tree as well as a $\bigvee $-tree.
\end{enumerate}
\end{punto}

The following results gives graphical sufficient/necessary conditions for a
given $X$-top lattice so that $X$ is completely normal.

\begin{thm}
\label{tcn}Let $\mathcal{L}=(L;\vee ,0;\wedge ,1)$ be an $X$-top lattice for
some $\emptyset \neq X\subseteq L\backslash \{1\}.$

\begin{enumerate}
\item If $X$ is a forest consisting of a finite number of pairwise \emph{%
strongly disjoint} $\bigwedge $-trees with finite base, then $X$ is
completely normal.

\item If $X$ is completely normal, then $X$ does \emph{not} contain any $%
\bigvee $-tree $\mathcal{V}$ with finite cover and $\left\vert Max(\mathcal{V%
})\right\vert \geq 2.$
\end{enumerate}
\end{thm}

\begin{Beweis}
\begin{enumerate}
\item Let $X$ be a forest consisting of a finite number of strongly disjoint
$\bigwedge $-trees with finite base $(\mathcal{T}_{\lambda _{1}},\mathfrak{m}%
_{i}),\cdots ,(\mathcal{T}_{\lambda _{n}},\mathfrak{m}_{n})$ and $%
I:=\{1,\cdots ,n\}.$

\textbf{Step 1:} $X$ is normal.

Let $V(a)$ and $V(b)$ be two \emph{disjoint} closed sets in $X.$ Then we can
find two \emph{disjoint} index sets $K,J\subseteq I$ such that $%
V(a)\subseteq \bigcup\limits_{i\in K}\mathcal{T}_{\lambda _{i}}$ and $%
V(b)\subseteq \bigcup\limits_{i\in J}\mathcal{T}_{\lambda _{i}}$ (if not,
then there exists some $1\leq j\leq n$ such that $\mathcal{T}_{\lambda
_{j}}\cap V(a)\neq \emptyset $ and $\mathcal{T}_{\lambda _{j}}\cap V(b)\neq
\emptyset ,$ whence $\mathfrak{m}_{j}\in V(a)\cap V(b)$, contradicting the
assumption that these sets are disjoint.

Since $\Lambda $ is \emph{finite}, $U:=\bigcap\limits_{i\in I\backslash
K}D(\bigwedge \mathcal{T}_{\lambda _{i}})$ and $V:=\bigcap\limits_{i\in
\Lambda \backslash J}D(\bigwedge \mathcal{T}_{\lambda _{i}})$ are two open
sets in $X.$

\textbf{Claim:} $V(a)\subseteq U,$ $V(b)\subseteq V$ and $U\cap V=\emptyset
. $

Let $x\in V(a).$ Since $V(a)\subseteq \bigcup\limits_{i\in K}\mathcal{T}%
_{\lambda _{i}},$ there exists $i\in K$ such that $x\in \mathcal{T}_{\lambda
_{i}}.$ If $x\notin U,$ then there exists $j\in I\backslash K$ such that $%
x\in V(\bigwedge \mathcal{T}_{\lambda _{j}}).$ So, $x\in V(\bigwedge
\mathcal{T}_{\lambda _{i}})\cap V(\bigwedge \mathcal{T}_{\lambda
_{j}})=\emptyset $ (since $i\neq j$), a contradiction. Similarly, we have $%
V(b)\subseteq V.$ Moreover, we have%
\begin{equation*}
U\cap V=\bigcap\limits_{i\in I\backslash (K\cap J)}D(\bigwedge \mathcal{T}%
_{\lambda _{i}})=\bigcap\limits_{i\in I}D(\bigwedge \mathcal{T}_{\lambda
_{i}})=X\backslash \bigcup\limits_{i\in I}V(\bigwedge \mathcal{T}_{\lambda
_{i}})=\emptyset .
\end{equation*}%
We conclude that, $X$ is normal.

\textbf{Step 2:} Every subspace of $X$ is normal.

Observe that removing any \emph{non maximal element} from any $\bigwedge $%
-tree with a finite base results in a (smaller) tree with finite base.
Moreover, removing the maximal element of a tree with finite base produces
two new disjoint $\bigwedge $-trees, each of finite base. Therefore, any
non-empty subspace $Y\subseteq X$ is again of the same type, whence normal
by Step 1.

\item Let $X$ be completely normal. Suppose that $X$ contains a dual tree $(%
\mathcal{V},m)$ with finite cover and $\left\vert Max(\mathcal{V}%
)\right\vert \geq 2.$ Let $\{y,z\}\subseteq Max(\mathcal{V})\subseteq X$ be
such that $m\lvertneqq y$ and $m\lvertneqq z.$ Notice that $Y:=\{m,y,z\}$ is
coatomic. Since $Y$ does \emph{not} have the pm-property, it follows by
Proposition \ref{txnpm} (1) that $Y$ is \emph{not} normal, a contradiction
to our assumption that $X$ is completely normal.$\blacksquare $
\end{enumerate}
\end{Beweis}

\begin{ex}
\label{not-pm}Consider the (bounded distributive)\ lattice $\mathcal{L}$ on $%
L=\{0,u,v,t,x,y,1\}$ given by the transitive closure of the relation
represented by the following diagram

\begin{figure}[]
\centering
\begin{equation*}
\xymatrix{ & & & { 1 \bullet} \ar@{-}[dll] \ar@{-}[drr] & & & \\ & x
{\bullet} \ar@{-}[drr] & & & & { \bullet} y \ar@{-}[dll] & \\ & & & {
\bullet} t \ar@{-}[dll] \ar@{-}[drr] & & & \\ & u { \bullet} \ar@{-}[drr] &
& & & v { \bullet} \ar@{-}[dll] & \\ & & & {\bullet} 0 & & & }
\end{equation*}%
\caption{An $X$-lattice with $X$ extremely non-normal but not anti-normal}
\label{fig.EN-not-AN}
\end{figure}
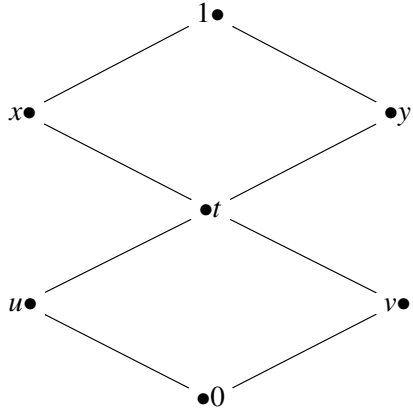
Let $X=\{x,u,y,v\}.$ Then $C^{X}(\mathcal{L})=\{0,t,x,y,u,v\}$ and $%
Max(X)=\{x,y\}.$ Clearly, $X=SI^{C^{X}(\mathcal{L})}(X),$ whence $L$ is an $%
X $-top lattice by Theorem \ref{xct}. The collection of closed sets is%
\begin{equation*}
\begin{tabular}{lllllllllll}
$V_{X}(0)$ & $=$ & $X,$ &  & $V_{X}(x)$ & $=$ & $\{x\},$ &  & $V_{X}(u)$ & $%
= $ & $\{u,x,y\},$ \\
$V_{X}(1)$ & $=$ & $\emptyset ,$ &  & $V_{X}(y)$ & $=$ & $\{y\},$ &  & $%
V_{X}(v)$ & $=$ & $\{v,y,x\},$ \\
$V_{X}(t)$ & $=$ & $\{x,y\}.$ &  &  &  &  &  &  &  &
\end{tabular}%
\end{equation*}%
However, $X$ does \emph{not} have the pm-property as $Max(u;X)=\{x,y\}$ ($%
=Max(v;X)$). By Corollary \ref{finite-pm-retract}, $X$ is \emph{not} normal.
We can double check this. The collection $\mathcal{O}(X)\mathcal{\ }$of open
sets in $X$ is%
\begin{equation*}
\begin{tabular}{lllllllllll}
$D_{X}(0)$ & $=$ & $\emptyset ,$ &  & $D_{X}(x)$ & $=$ & $\{y,u,v\},$ &  & $%
D_{X}(u)$ & $=$ & $\{v\},$ \\
$D_{X}(1)$ & $=$ & $X,$ &  & $D_{X}(y)$ & $=$ & $\{x,u,v\},$ &  & $D_{X}(v)$
& $=$ & $\{u\},$ \\
$D_{X}(t)$ & $=$ & $\{u,v\}.$ &  &  &  &  &  &  &  &
\end{tabular}%
\end{equation*}%
Notice that $Cl(X)=\{x,y\}=Max(X)$ and $Iso(X)=\{u,v\}.$ So, $X=Cl(X)\cup
Iso(X),$ i.e., $X$ is $T_{\frac{1}{2}}.$ However, $X$ is not $T_{\frac{3}{4}%
},$ e.g., $x$ is \emph{not} regular open as $int(\overline{\{x\}}%
)=int(\{x\})=\emptyset $ (cf. \cite[Theorem 3.20 (2)]{AA}). Notice that $X$
is \emph{not} normal as the closed sets $\{x\},$ $\{y\}$ cannot be separated
by disjoint open sets. We demonstrate also that $X$ is \emph{not}
max-retractable. Suppose that $f:X\longrightarrow Max(X)$ were a retraction.
Since $f_{|_{Max(X)}}=id_{X},$ we would have $f(x)=x$ and $f(y)=y.$ Since $%
y\in V(u),$ it would follow then that $f(y)\in f(V(u))\overset{\text{Remark %
\ref{ret-unique}}}{=}f(x)=x,$ a contradiction.

Moreover, this is an example of an extremely non-normal space which is not
anti-normal. Notice that $\mathcal{T}(\mathbf{N};X)=\{(\{x\}\times \{y\})\}$
while $\mathcal{S}(\mathbf{N};X)=\emptyset ,$ i.e. $X$ is extremely
non-normal. On the other hand, the subspaces $Y=\{u,v,y\}$ and $Z=\{u,v,x\}$
are obviously normal, whence $X$ is not anti-normal.$\blacksquare $
\end{ex}

\begin{rem}
Example \ref{not-pm} illustrates also that the assumption that the trees in
the forest in Theorem \ref{tcn} (1)\ are \emph{strongly disjoint} cannot be
even weakened by replacing it with the assumption that the trees are just
\emph{disjoint}. One might think that $X=\mathcal{T}_{2}\cup \mathcal{T}%
_{2}^{\prime }=\{x,u,v\}\bigcup \{y,u,v\}$ would be completely normal if $%
\mathcal{T}_{2}$ and $\mathcal{T}_{2}^{\prime }$ were disjoint. But that is
not the case as $X=\mathcal{T}_{1}\bigsqcup \mathcal{T}_{1}^{\prime
}=\{u,x\}\bigsqcup \{v,y\},$ a \emph{disjoint} union of trees with finite
bases. However, $\mathcal{T}_{1}$ and $\mathcal{T}_{1}^{\prime }$ are \emph{%
not strongly disjoint} as
\begin{equation*}
V(\bigwedge \mathcal{T}_{1})\cap V(\bigwedge \mathcal{T}_{1}^{\prime
})=V(u)\cap V(v)=\{x,y\}\neq \emptyset .\blacksquare
\end{equation*}
\end{rem}

The following result provides a necessary condition for an $X$-top lattice
so that $X$ is \emph{perfectly normal}.

\begin{prop}
\label{ppn}Let $\mathcal{L}=(L;\vee ,0;\wedge ,1)$ be an $X$-top lattice for
some $\emptyset \neq X\subseteq L\backslash \{1\}.$

\begin{enumerate}
\item If $X$ is perfectly normal, then $X$ does not contain any $\mathcal{C}%
_{2}$ as a subset, i.e., $K.\dim (X)=0$.

\item $X$ is $T_{6}$ if and only if $X$ is perfectly normal.
\end{enumerate}
\end{prop}

\begin{Beweis}
\begin{enumerate}
\item Let $X$ be perfectly normal. Suppose that $X$ contains a $\mathcal{C}%
_{2}$ (i.e., there exist $x,y\in X$ such that $x\lvertneqq y$). Notice that
if $D(a)$ is any open subset of $X$ containing $V(y),$ then $x\in D(a).$ It
follows that the closed set $V(y)$ cannot be expressed as an intersection of
open sets, contradicting the assumption that $X$ is perfectly normal.

\item Let $X$ be perfectly normal. By (1), $X$ does not contain a $\mathcal{C%
}_{2},$ whence $K.\dim (X)=0$ or equivalently $X$ is $T_{1}$ by Proposition %
\ref{dim-0} (4). Notice that this follows also from the fact that a $T_{0}$
perfectly normal space is $T_{1}$ (cf. Remark \ref{rem-T3}).$\blacksquare $
\end{enumerate}
\end{Beweis}

\begin{lem}
\label{localization-AM}\emph{(cf. \cite[Exercise 1.22 $\&\ $Proposition 3.11]%
{AM1969})}\ Let $R$ be a commutative ring.

\begin{enumerate}
\item Let $R=\prod\limits_{i=1}^{n}R_{i}$ (a finite direct product of
rings). Then $Spec(R)\approx \bigsqcup\limits_{i=1}^{n}Spec(R_{i})$ (i.e., $%
Spec(R)$ is homeomorphic to the disjoint union of the prime spectra of the
rings $R_{1},\cdots ,R_{n}$).

\item If $S\subseteq R$ is a multiplicatively closed set, then the is a 1-1
correspondence%
\begin{equation*}
\{P\in Spec(R)\text{ }|\text{ }P\cap S=\emptyset \}\longleftrightarrow
Spec(S^{-1}R),\text{ }P\longmapsto S^{-1}P.
\end{equation*}
\end{enumerate}
\end{lem}

\begin{ex}
\label{stongly-disjoint}Consider the (bounded distributive)\ lattice $%
\mathcal{L}$:
\begin{figure}[]
\centering
\begin{equation*}
\xymatrix{ & & & { 1 \bullet} \ar@{-}[dll] \ar@{-}[drr] & & & \\ & x
{\bullet} \ar@{-}[ddl] \ar@{-}[drr] & & & & { \bullet} y \ar@{-}[ddr]
\ar@{-}[dll] & \\ & & & { \bullet} t \ar@{-}[ddl] \ar@{-}[ddr] & & & \\ u {
\bullet} \ar@{-}[drr] & & & & & & v { \bullet} \ar@{-}[dll] \\ & & {\bullet}
z & & {\bullet} w & & \\ & & & {\bullet} 0 \ar@{-}[ul] \ar@{-}[ur] & & & }
\end{equation*}%
\caption{An $X$-top lattice with $X$ completely normal but not perfectly
normal}
\label{X-cn-not-pn}
\end{figure}
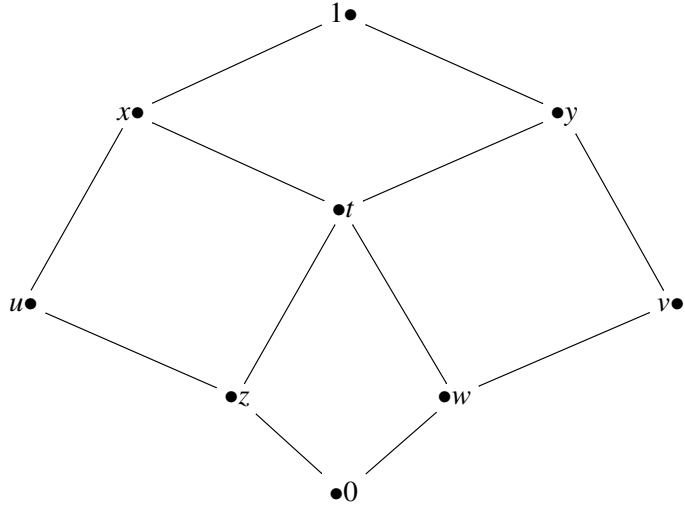
Let $X=\{x,u,y,v\}.$ Then $C^{X}(\mathcal{L})=\{0,t,x,y,u,v,w,z\}$ and $%
Max(X)=\{x,y\}.$ Clearly, $X=SI^{C^{X}(\mathcal{L})}(X),$ whence $L$ is an $%
X $-top lattice by Theorem \ref{xct}. The collection of closed sets is%
\begin{equation*}
\begin{tabular}{lllllllllll}
$V_{X}(0)$ & $=$ & $X,$ &  & $V_{X}(x)$ & $=$ & $\{x\},$ &  & $V_{X}(v)$ & $%
= $ & $\{v,y\},$ \\
$V_{X}(1)$ & $=$ & $\emptyset ,$ &  & $V_{X}(y)$ & $=$ & $\{y\},$ &  & $%
V_{X}(z)$ & $=$ & $\{u,x,y\},$ \\
$V_{X}(t)$ & $=$ & $\{x,y\}.$ &  & $V_{X}(u)$ & $=$ & $\{u,x\},$ &  & $%
V_{X}(w)$ & $=$ & $\{v,x,y\}.$%
\end{tabular}%
\end{equation*}%
Notice that $X=\{u,x\}\bigsqcup \{u,v\}=\mathcal{T}_{1}\bigsqcup \mathcal{T}%
_{1}^{\prime },$ a disjoint union of two \emph{strongly disjoint} trees ($%
V_{X}(\bigvee \mathcal{T}_{1})\cap V_{X}(\bigvee \mathcal{T}_{1}^{\prime
})=V_{X}(u)\cap V_{X}(v)=\{u,x\}\cap \{v,y\}=\emptyset $). Moreover, $\dim
(X)=1.$ It follows, by Theorem \ref{tcn} (1)\ and Proposition \ref{ppn} (1),
that $X$ is completely normal but \emph{not} perfectly normal. The
collection of open sets of $X$ is
\begin{equation*}
\begin{tabular}{lllllllllll}
$D_{X}(0)$ & $=$ & $X,$ &  & $D_{X}(x)$ & $=$ & $\{y,u,v\},$ &  & $D_{X}(v)$
& $=$ & $\{u,x\},$ \\
$D_{X}(1)$ & $=$ & $\emptyset ,$ &  & $D_{X}(y)$ & $=$ & $\{x,u,v\},$ &  & $%
D_{X}(z)$ & $=$ & $\{v\},$ \\
$D_{X}(t)$ & $=$ & $\{u,v\},$ &  & $D_{X}(u)$ & $=$ & $\{v,y\},$ &  & $%
D_{X}(w)$ & $=$ & $\{u\}.$%
\end{tabular}%
\end{equation*}%
Notice that $Cl(X)=\{x,y\}=Max(X)$ and $Iso(X)=\{u,v\}.$ So, $X=Cl(X)\cup
Iso(X),$ i.e., $X$ is $T_{\frac{1}{2}}.$ However, $X$ is not $T_{\frac{3}{4}%
},$ e.g., $x$ is \emph{not} regular open as $int(\overline{\{x\}}%
)=int(\{x\})=\emptyset $ (cf. \cite[Theorem 3.20 (2)]{AA}). One can double
check that $X$ is completely normal but not a $G_{\delta }$-space (e.g., all
open sets containing $x$ contain $u$ as well, whence the closed set $\{x\}$
is not a $G_{\delta }$-set). It follows by Lemma \ref{Vedenissoff} that $X$
is not perfectly normal. Notice also that $X$ is not regular:\ $%
(\{x,y\},u)\in \mathcal{T}(\mathbf{R};X);$ however, $\{x,y\}$ and $u$ cannot
be separated by disjoint open sets.

Notice that $X$ is spectral (being finite and $T_{0}$). By Lemma \ref%
{localization-AM}, $X\approx Spec(R),$ where $R$ can be chosen so that $%
R=D_{1}\times D_{2}$ and $(D_{1},\mathfrak{m}_{1}),$ $(D_{2},\mathfrak{m}%
_{2})$ are two DVRs. In this case, we have%
\begin{equation*}
X\approx Spec(D_{1})\bigsqcup Spec(D_{2})=\{0_{1}\times D_{2},\mathfrak{m}%
_{1}\times D_{2},D_{1}\times 0_{2},D_{1}\times \mathfrak{m}_{2}\}.
\end{equation*}%
Moreover, making use of Proposition \ref{Spec(B(n,i))} (3), $X\approx
Spec(S),$ where $S=S_{1}\times S_{2}$ and $S_{i}=(B(n_{i},n_{i}-1),\mathfrak{%
m}_{n_{i}})$ with $n_{1},n_{2}\geq 3$ is the local Alarcon-Anderson
semidomain \cite{AA1994} with maximal ideal $\mathfrak{m}_{n_{i}}:=\{0,2,%
\cdots ,n_{i}-1\}.$ In this case%
\begin{equation*}
X\approx Spec(S_{1})\bigsqcup Spec(S_{2})=\{0_{1}\times S_{2},\mathfrak{m}%
_{n_{1}}\times S_{2},S_{1}\times 0_{2},S_{1}\times \mathfrak{m}_{n_{2}}\}.
\end{equation*}
\end{ex}

\section{Regular $X$-top lattices}

In this section, we study \emph{regular }and\emph{\ completely regular }$X$%
-top lattices.

\bigskip

Since any $X$-top lattice is $T_{0}$ by Proposition \ref{dim-0}, we have the
following result by Remark \ref{rem-T3} (3). However, we include a short
proof for the sake of completeness.

\begin{lem}
\label{rcr}Let $\mathcal{L}=(L;\vee ,0;\wedge ,1)$ be a $X$-top lattice.

\begin{enumerate}
\item $X$ is regular iff $X$ is $T_{3};$

\item $X$ is completely regular if and only if $X$ is $T_{3\frac{1}{2}}.$
\end{enumerate}
\end{lem}

\begin{Beweis}
Let $X$ be regular. It is enough to show that $K.\dim (X)=0.$ Suppose that
there exists $x\lvertneqq y$ in $X$. Then $V(y)$ is a closed set and $%
x\notin V(y)$. However, any open set $D(a)\ $containing $V(y)$ must contain $%
x$ (a contradiction to the \emph{regularity} of $X$).$\blacksquare $
\end{Beweis}

\begin{prop}
\label{cpmpact-T3-T4}Let $\mathcal{L}=(L;\vee ,0;\wedge ,1)$ be an $X$-top
lattice for some $\emptyset \neq X\subseteq L\backslash \{1\}.$ If $X$ is
\emph{compact}, then the following are equivalent:

\begin{enumerate}
\item $X$ is regular;

\item $X$ is $T_{3};$

\item $X$ is $T_{2\frac{1}{2}};$

\item $X$ is $T_{2};$

\item $X$ is $T_{4};$

\item $X$ is $T_{3\frac{1}{2}};$

\item $X$ is $T_{1}$ and quasi-Hausdorff;

\item $K.\dim (X)=0$ and quasi-Hausdorff.

If moreover, $X$ is coatomic, then these are equivalent to

\item $X$ is normal and Jacobson.
\end{enumerate}
\end{prop}

\begin{Beweis}
$(1\Longleftrightarrow 2)$ This is Lemma \ref{rcr} (1).

$(2\Longrightarrow 3\Longrightarrow 4)\ $and $(5\Longrightarrow
6\Longrightarrow 1)$ Obvious.

$(4\Longrightarrow 5)$ Any compact Hausdorff space is normal (e.g., \cite[%
Theorem 17.10]{Wil1970}).

$(4\Longleftrightarrow 7)$ This is Lemma \ref{T1-qH-T2} (3).

$(4\Longleftrightarrow 8)$ This is Proposition \ref{dim-0} (5).

$(5\Longrightarrow 9)$ Obvious.

Assume that $X$ is \emph{coatomic}.

$(9\Longrightarrow 5)$ Since $X$ coatomic and normal, it follows by
Proposition \ref{txnpm} that $X$ has the pm-property. Since $X$ is
(moreover) Jacobson, it follows that $X=Max(X),$ whence $X$ is $T_{1}$ by
Proposition \ref{dim-0} (4).$\blacksquare $
\end{Beweis}

\begin{defn}
(cf. \cite[29.4]{Wil1970})\ We say that a topological space $X$ is \textbf{%
inductively zero-dimensional }iff $X$ has a base of \emph{clopen sets}
(equivalently, $ind.\dim (X)=0,$ where $ind.\dim (X)$ is the so-called
\textbf{small inductive dimension} of $X$ \cite[page 105]{AP1990}).
\end{defn}

\begin{defn}
(cf. \cite{Joh1982}) A \textbf{Stone space (Boolean space, profinite space})
is a topological space $X$ that satisfies any, hence all, of the following
equivalent conditions:

\begin{enumerate}
\item $X$ is $T_{0},$ inductively zero-dimensional space and compact;

\item $X$ is homeomorphic to a projective limit of \emph{finite discrete}
spaces.
\end{enumerate}
\end{defn}

Combining Proposition \ref{cpmpact-T3-T4} and \cite[Propostion 2.15]{AA}, we
obtain further characterizations of $X$-top lattices for which $X$ is a
Stone space.

\begin{thm}
\label{abs-flat-regular}Let $\mathcal{L}=(L;\vee ,0;\wedge ,1)$ be an $X$%
-top lattice for some $\emptyset \neq X\subseteq L\backslash \{1\}.$ The
following are equivalent:

\begin{enumerate}
\item $X$ is a Stone space;

\item $X$ is spectral and $T_{2};$

\item $X$ is spectral and regular;

\item $X$ is spectral and $T_{3};$

\item $X$ is spectral and $T_{4};$

\item $X$ is spectral and $T_{3\frac{1}{2}};$

\item $X$ is spectral, Jacobson and normal;

\item $X$ is spectral, Jacobson and has the pm-property;

\item $X$ is spectral and $K.\dim (X)=0;$

\item $X$ is spectral, dual Jacobson and has the m-property;

\item $X$ is spectral and $K.\dim (X)=0;$

\item $X$ is homeomorphic to $Spec(R)$ for some commutative Jacobson
pm-(semi)ring.

\item $X$ is homeomorphic to $Spec(R)$ for some commutative dual Jacobson
m-(semi)ring.
\end{enumerate}
\end{thm}

\begin{defn}
We call a (semi)ring $R$

\textbf{von Neumann regular} iff for every $a\in R$ there exists $b$ such
that $a=aba;$

$\pi $\textbf{-regular} iff for every $a\in R$ there exists $b\in R$ and $%
n\geq 1$ such that $a^{n}=a^{n}ba^{n};$

\textbf{reduced }iff $Nil(R):=\{a\in R$ $\mid $ $a^{n}=0$ for some $n\geq
1\}=0.$
\end{defn}

Combining Theorem \ref{abs-flat-regular} with well-known characterizations
of Krull $0$-dimensional commutative rings (e.g., \cite{Sto1968}, \cite[%
Exercise 3.11]{AM1969}, \cite{Gil2000}), we obtain:

\begin{cor}
\label{vN-regular}The following are equivalent for a commutative ring $R:$

\begin{enumerate}
\item $R$ is von Neumann regular;

\item $R$ is reduced and $\pi $-regular;

\item $R$ is reduced and $Spec(R)$ is a Stone space;

\item $R$ is reduced and $Spec(R)$ is $T_{2};$

\item $R$ is reduced and $Spec(R)$ is $T_{1};$

\item $R$ is reduced and $K.\dim (R)=0;$

\item $R$ is a reduced Jacobson pm-ring;

\item $R$ is a reduced dual Jacobson m-ring;

\item $R$ is reduced and $Spec(R)$ is regular ($T_{3}$);

\item $R$ is a reduced Jacobson ring and $Spec(R)$ is normal ($T_{4}$).
\end{enumerate}
\end{cor}

\bigskip

The following example illustrates that Corollary \ref{vN-regular} does not
apply for \emph{proper} commutative semirings (that are \emph{not} rings):

\bigskip

\begin{ex}
\label{B(3,1)}Consider the Alarcon-Anderson semidomain $B(3,1)=\{0,1,2\}.$
It can be easily seen that $B(3,1)$ is von Neumann regular (notice that $%
2^{2}\cdot 1=2$ in $B(3,1)$). Notice that $Spec(B(3,1))=\{\{0\},\{0,2\}\}$
(cf. \ref{Spec(B(n,i))}), is completely normal by Theorem \ref{tcn} but
\emph{not} perfectly normal as it is a $\mathcal{C}_{2}$ (cf. Proposition %
\ref{ppn}). Clearly, $B(3,1)$ is reduced and normal but \emph{not} a
Jacobson semiring. Being finite, $B(3,1)$ is an Artinian (and Noetherian)
semiring, but $K.\dim (Spec(B(3,1))\neq 0.$ It follows by Lemma \ref{rcr}
that $Spec(B(3,1))$ is \emph{not} regular (priori \emph{not} completely
regular). Notice that $Spec(B(3,1))\approx G$ in Example \ref{ex-DVR} below
an is homeomorphic to the Sierpi\'{n}ski space (the smallest example of a
topological space which is neither discrete nor indiscrete).$\blacksquare $
\end{ex}

\subsection*{Examples and Counterexamples}

\bigskip

We devote the rest of this article to examples and counterexamples.

\bigskip

\begin{notatation}
We denote by $\mathbb{P}=\{2,3,5,\cdots \}$ the set of prime numbers. For a
positive integer $m,$ we denote by $\mathbb{P}(m)$ the set of prime divisors
of $m$ and set $\omega (m):=\left\vert \mathbb{P}(m)\right\vert .$ For a
prime $p\in \mathbb{P}$, we denote the by $\mathbb{F}_{(p)}$ the
localization of $\mathbb{Z}$ at the prime ideal $p\mathbb{Z}.$ We denote
with $\mathcal{D}_{n+2}$ the \emph{generalized diamond}:\ a poset with a
\emph{largest element} $\mathfrak{m}$ and a \emph{smallest element} $m$ in
addition to $n$ (intermediate)\ incomparable elements which we denote by $%
\mathcal{P}_{n}.$
\end{notatation}

\bigskip

Fix the \emph{distributive} lattice $\mathcal{L}=(L;\vee ,0;\wedge ,1)$
given by the transitive completion of the partial order on $%
L=\{0,z,x,y,w,1\} $ represented in the diagram below

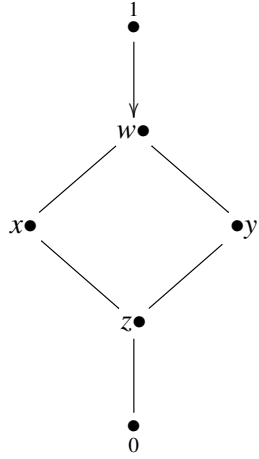
\begin{figure}[]
\centering
\begin{equation*}
\xymatrix{ & \overset{1}{\bullet} \ar[d] & \\ & w { \bullet} \ar@{-}[dl]
\ar@{-}[dr] & \\ x { \bullet} \ar@{-}[dr] & & { \bullet} y \ar@{-}[dl] \\ &
z { \bullet} \ar@{-}[d] & \\ & \underset{0}{\bullet} & }
\end{equation*}%
\caption{A distributive modular lattice}
\label{fig:LN}
\end{figure}

We will investigate five subsets of $L\backslash \{1\}$ for each of which
the poset $(L,\leq )$ induces a Zariski-like topology.

\bigskip

While ultraconnected spaces are \emph{trivially} normal, the following
example illustrates that such spaces are \emph{not necessarily} completely
normal (priori \emph{not necessarily} perfectly normal). Moreover, it
provides an $X$-top lattice for which $X$ is a $T_{\frac{1}{2}}$ normal
space but \emph{not} regular ($X$ is even \emph{extremely non-regular}).

\begin{ex}
\label{normal-not-regular}Consider the lattice $\mathcal{L}$ in Figure \ref%
{fig:LN}. Set $X:=\{0,x,y,w\}$ (the dotted lines indicate that $0\lvertneqq
x\wedge y$ in $\mathcal{L}$):
\begin{figure}[]
\centering
\begin{equation*}
\xymatrix{ & w { \bullet} \ar@{-}[dl] \ar@{-}[dr] & \\ x { \bullet}
\ar@{.}[dr] & & { \bullet} y \ar@{.}[dl] \\ & \underset{0}{\bullet} & }
\end{equation*}%
\caption{$X$ is normal but extremely non-regular}
\label{n-not-cn-enr}
\end{figure}
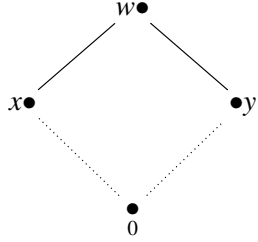

Notice that $C^{X}(\mathcal{L})=\{0,z,x,y,w\}$. It is clear that every
element $p\in X$ is strongly $C^{X}(\mathcal{L})$-irreducible. By Theorem %
\ref{xct}, $X$ attains a Zariski-like topology induced by the poset $(L,\leq
).$ One can double check:\ The collection $V_{X}(\mathcal{L})$ of closed
varieties consists of%
\begin{equation*}
\begin{tabular}{lllllllllll}
$V_{X}(0)$ & $=$ & $X,$ &  & $V_{X}(z)$ & $=$ & $\{x,y,w\},$ &  & $V_{X}(y)$
& $=$ & $\{y,w\},$ \\
$V_{X}(1)$ & $=$ & $\emptyset ,$ &  & $V_{X}(x)$ & $=$ & $\{x,w\},$ &  & $%
\text{ }V_{X}(w)$ & $=$ & $\{w\}.$%
\end{tabular}%
\end{equation*}%
It's obvious that $V_{X}(\mathcal{L})\ $is closed under finite unions.
Notice that $X$ is colocal with $Min(X)=\{0\}$ and local with $Max(X)=\{w\}$%
. Since $X$ is finite and has the pm-property, it follows by Corollary \ref%
{finite-pm-retract} that $X$ is normal. However, $X$ is \emph{not}
completely normal since $\{0,x,y\}\subseteq X$ is a $\mathcal{V}_{2}$ (cf.
Theorem \ref{tcn}).

Since $\dim (X)\neq 0,$ it follows by Proposition \ref{cpmpact-T3-T4} that $%
X $ is \emph{not} regular whence \emph{neither} completely regular \emph{nor}
perfectly normal (cf. Proposition \ref{ppn} (1)). In fact, $X$ is far away
from being regular:\ $X$ is hyperconnected (irreducible), whence extremely
non-$T_{2}$ by Lemma \ref{anti-hyper} (1). Since $\mathcal{T}(\mathbf{R}%
;X)\neq \emptyset ,$ it follows that $X$ is \emph{extremely non-regular} by
Lemma \ref{hyp-e-non-regular}. Moreover, $X$ is far away from being $T_{1}$
as $X$ is \emph{not even} $T_{\frac{1}{4}}$ (cf. Proposition \ref{dim-0}
(6)).

We can double check these observations: Notice that $X$ is ultraconnected
(as $w\in C$ for every $\emptyset \neq C\underset{\text{closed}}{\subseteq }%
X $), whence $\mathcal{T}(\mathbf{N};X)=\emptyset $ and consequently $X$ is
\emph{trivially} normal. To see that $X$ is \emph{not} completely normal,
notice that $A=\{x\}$ and $B=\{y\}$ are \emph{separated}:%
\begin{equation*}
A\cap \overline{B}=\{x\}\cap \{y,w\}=\emptyset =\{x,w\}\cap \{y\}=\overline{A%
}\cap B.
\end{equation*}%
However, $A$ and $B$ cannot be separated by open sets as the collection $%
\mathcal{O}(X)$ of open subsets of $X$ consists of%
\begin{equation*}
\begin{tabular}{lllllllllll}
$D_{X}(0)$ & $=$ & $\emptyset ,$ &  & $D_{X}(z)$ & $=$ & $\{0\},$ &  & $%
D_{X}(y)$ & $=$ & $\{0,x\},$ \\
$D_{X}(1)$ & $=$ & $X,$ &  & $D_{X}(x)$ & $=$ & $\{0,y\},$ &  & $\text{ }%
D_{X}(w)$ & $=$ & $\{0,x,y\}.$%
\end{tabular}%
\end{equation*}%
It follows, by Lemma \ref{CN-separated}, that $X$ is \emph{not} completely
normal. Moreover, $X$ is \emph{not} a $G_{\delta }$-space (\emph{e.g.}, the
closed set $\{w\}$ is not contained in any open set), whence $X$ is \emph{not%
} perfectly normal by Lemma \ref{Vedenissoff}.

Since $X$ is finite and $T_{0},$ it follows by Lemma \ref{dim-0} (2) that $X$
is spectral. In fact, $X\approx Spec(S),$ where we may choose (cf. \ref%
{Spec(B(n,i))}):%
\begin{eqnarray*}
S &=&B(n,i),\text{ }n\geq 8,\text{ }2\leq i\leq n-i\text{ and }\mathbb{P}%
(n-i)=\{p,q\}; \\
Spec(S) &=&\{0,pS,qS,\mathfrak{m}_{n}\},\text{ where }\mathfrak{m}%
_{n}=\{0,2,3,\cdots ,n-1\}\text{).}\blacksquare
\end{eqnarray*}
\end{ex}

\bigskip

In what follows, we give an example of a $Y$-top lattice for which $Y$ is a $%
T_{\frac{3}{4}}$ completely normal space but \emph{not} perfectly normal.

\bigskip

\begin{ex}
\label{CN-not-PN}Consider the lattice $\mathcal{L}$ in Figure \ref{fig:LN}.
Consider $Y=\{x,y,w\}:$
\begin{figure}[]
\centering
\begin{equation*}
\xymatrix{& w { \bullet} \ar@{-}[dl] \ar@{-}[dr] & \\ x { \bullet} & & {
\bullet} y }
\end{equation*}%
\caption{$Y$ is completely normal $\&$ extremely non-regular but not
extremely non-$T_2$}
\label{Y-cn}
\end{figure}
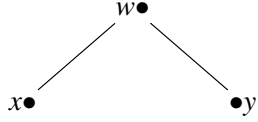
Notice that $C^{Y}(\mathcal{L})=\{z,x,y,w\}.$ It follows, by Example \ref%
{normal-not-regular} and Corollary \ref{Y-top}, that $Y$ attains a
Zariski-like topology induced by the poset $(L,\leq ).$ One can double
check:\ the proposed collection of closed sets consists of%
\begin{equation*}
\begin{tabular}{lllllllllll}
$V_{Y}(0)$ & $=$ & $Y,$ &  & $V_{Y}(z)$ & $=$ & $Y,$ &  & $V_{Y}(y)$ & $=$ &
$\{y,w\},$ \\
$V_{Y}(1)$ & $=$ & $\emptyset ,$ &  & $V_{Y}(x)$ & $=$ & $\{x,w\},$ &  & $%
\text{ }V_{Y}(w)$ & $=$ & $\{w\}.$%
\end{tabular}%
\end{equation*}%
It's clear that $V_{Y}(\mathcal{L})\ $is closed under finite unions. Notice
that, $Y$ is $\mathcal{T}_{2}$ and $K.\dim (Y)=1.$ Hence, $Y$ is completely
normal but \emph{not} perfectly normal (by Theorem \ref{tcn} and Proposition %
\ref{ppn}). By Lemma \ref{rcr}, $Y$ is \emph{not} regular as $K.\dim (Y)\neq
0$ (priori \emph{neither} completely regular \emph{nor} perfectly normal).
Notice that $Y$ is $T_{\frac{3}{4}}$ but not $T_{1}$ by \cite[Theorem 3.20
(1)]{AA}.

The collection $\mathcal{O}(Y)$ of open subsets of $Y$ consists of%
\begin{equation*}
\begin{tabular}{lllllllllll}
$D_{Y}(0)$ & $=$ & $\emptyset ,$ &  & $D_{Y}(z)$ & $=$ & $\emptyset ,$ &  & $%
D_{Y}(y)$ & $=$ & $\{x\},$ \\
$D_{Y}(1)$ & $=$ & $Y,$ &  & $D_{Y}(x)$ & $=$ & $\{y\},$ &  & $\text{ }%
D_{Y}(w)$ & $=$ & $\{x,y\}.$%
\end{tabular}%
\end{equation*}%
Since $\mathcal{T}(\mathbf{R};Y)=\{(\{x,w\},y),$ $(\{y,w\},x),$ $(\{w\},x),$
$(\{w\},y)\},$ one can easily see that $\mathcal{S}(\mathbf{R};Y)=\emptyset
, $ i.e. $Y$ is \emph{even} extremely non-regular. Notice that $X$ is \emph{%
not} extremely non-Hausdorff, as $X$ is \emph{not} hyperconnected (i.e. the
converse of Lemma \ref{hyp-e-non-regular} (2-a) is \emph{not} true).

Since $Y$ is finite and $T_{0}$, it follows by Lemma \ref{dim-0} (2) that $Y$
is spectral. In fact $X\approx Spec(R),$ where we can choose (for any field $%
\mathbb{F}$):
\begin{equation*}
R=\mathbb{F}[[x,y]]/(xy)\text{ with }Spec(R)=\{(x),(y),(x,y)\}.\blacksquare
\end{equation*}
\end{ex}

\bigskip

In what follows, we give an example of a $Q$-top lattice for which $Q$ does
\emph{not} have the pm-property, whence $Q$ is \emph{not} normal (even
extremely non-normal). Moreover, $Q$ is the unique anti-normal space (up to
homeomorphism).

\bigskip

\begin{ex}
\label{anti-normal}Consider the lattice $\mathcal{L}$ in Figure \ref{fig:LN}%
. Set $Q:=\{0,x,y\}$ (the dotted lines indicate that $0\lvertneqq x\wedge y$
in $\mathcal{L}$):
\begin{figure}[]
\centering
\begin{equation*}
\xymatrix{ x { \bullet} \ar@{.}[dr] & & { \bullet} y \ar@{.}[dl] \\ &
\underset{0}{\bullet} & }
\end{equation*}%
\caption{$Q$ is the unique anti-normal space up to homeomorphism}
\label{Q-anti-normal}
\end{figure}
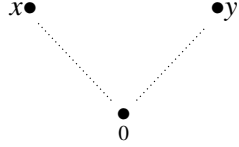

Notice that $C^{Q}(\mathcal{L})=\{0,z,x,y\}.$ It follows, by Example \ref%
{normal-not-regular} and Corollary \ref{Y-top}, that $Q$ attains a
Zariski-like topology induced by the poset $(L,\leq ).$ One can double
check:\ the proposed collection of closed sets consists of%
\begin{equation*}
\begin{tabular}{lllllllllll}
$V_{Q}(0)$ & $=$ & $Q,$ &  & $V_{Q}(z)$ & $=$ & $\{x,y\},$ &  & $V_{Q}(y)$ &
$=$ & $\{y\},$ \\
$V_{Q}(1)$ & $=$ & $\emptyset ,$ &  & $V_{Q}(x)$ & $=$ & $\{x\},$ &  & $%
\text{ }V_{Q}(w)$ & $=$ & $\emptyset .$%
\end{tabular}%
\end{equation*}%
It's clear that $V_{Y}(\mathcal{L})\ $is closed under finite unions. Notice
that, $Q$ is $\mathcal{V}_{2}$ and $K.\dim (Q)=1.$ Notice that $Q$ is
colocal with $Min(X)=\{0\}$ but \emph{not} local as $Max(Q)=\{x,y\}.$
Moreover, $Q$ is coatomic but does \emph{not} have the pm-property since $%
\left\vert Max(0;Q)\right\vert >1.$ It follows, by Proposition \ref{txnpm}
(1), that $Q$ is \emph{not} normal (priori \emph{neither} completely normal
\emph{nor} perfectly normal). One can also see that $X$ is \emph{not }%
completely normal since $X$ is $\mathcal{V}_{2}$ (cf. Theorem \ref{tcn}
(2)), and that $Q$ is \emph{not} perfectly normal as $Q$ contains a $%
\mathcal{C}_{2}$ (cf. Proposition \ref{ppn} (1)).

Moreover, it follows by Proposition \ref{rcr}, that $Y$ is \emph{not}
regular as $K.\dim (Q)\neq 0$. The collection $\mathcal{O}(Q)$ of open
subsets of $Y$ consists of%
\begin{equation*}
\begin{tabular}{lllllllllll}
$D_{Q}(0)$ & $=$ & $Q,$ &  & $D_{Q}(z)$ & $=$ & $\{0\},$ &  & $D_{Q}(y)$ & $%
= $ & $\{0,x\},$ \\
$D_{Q}(1)$ & $=$ & $\emptyset ,$ &  & $D_{Q}(x)$ & $=$ & $\{0,y\},$ &  & $%
\text{ }D_{Q}(w)$ & $=$ & $Q.$%
\end{tabular}%
\end{equation*}%
Clearly, $Q$ is hyperconnected (irreducible), whence $Q$ is extremely
non-Hausdorff. Since $\mathcal{T}(\mathbf{R};Q)\neq \emptyset $ and $%
\mathcal{T}(\mathbf{N};Q)\neq \emptyset ,$ it follows by Lemma \ref%
{anti-hyper} (2) that $Q$ is extremely non-regular and extremely non-normal.
In fact, $Q$ is the \emph{unique} anti-normal topological spaces (up to
homeomorphism) as show in \cite{GRV1981}. Notice that $Q$ is, as well, far
away from being $T_{1}$ as $Q$ is $T_{\frac{1}{2}}$ but \emph{not} $T_{\frac{%
3}{4}}$ by \cite[Corollary 3.11. $\&\ $Theorem 3.20 (2)]{AA}.

Since $Q$ is finite and $T_{0}$, it follows by Lemma \ref{dim-0} (2) that $Q$
is spectral:

\begin{itemize}
\item $Q\approx Spec(R),$ where
\begin{eqnarray*}
R &=&\mathbb{Z}_{(2)}\cap \mathbb{Z}_{(3)}=\{\frac{m}{n}\in \mathbb{Q}%
:g.c.d.(n,6)=1\}; \\
Spec(R) &=&\{0,\text{ }2R,\text{ }3R\}.
\end{eqnarray*}

\item $Q\approx Spec(S),$ where we can choose $S$ to be any of the
Alarcon-Anderson semidomains%
\begin{eqnarray*}
S &=&B(n,1)\text{ with }n\geq 7,\text{ }\omega (n-1)=2\text{ and }\mathbb{P}%
(n-1)=\{p,q\}; \\
Spec(S) &=&\{0,pS,qS\}.\blacksquare
\end{eqnarray*}
\end{itemize}
\end{ex}

\bigskip

\begin{ex}
\label{orerfectly normal}Consider the lattice $\mathcal{L}$ in Figure \ref%
{fig:LN}. Set $H:=\{x,y\}$:
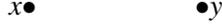
\begin{figure}[]
\centering
\begin{equation*}
\xymatrix{ x { \bullet} & & { \bullet} y}
\end{equation*}%
\caption{A perfectly normal ($T_6$) space}
\label{T6}
\end{figure}
Notice that $C^{H}(\mathcal{L})=\{z,x,y\}$. It follows, by Example \ref%
{normal-not-regular} and Corollary \ref{Y-top}, that $H$ attains a
Zariski-like topology induced by the poset $(L,\leq ).$ One can double
check:\ The collection of open sets of $H$ is given by%
\begin{equation*}
\begin{tabular}{lllllllllll}
$D_{H}(0)$ & $=$ & $\emptyset ,$ &  & $D_{H}(z)$ & $=$ & $\emptyset ,$ &  & $%
D_{H}(y)$ & $=$ & $\{x\},$ \\
$D_{H}(1)$ & $=$ & $H,$ &  & $D_{H}(x)$ & $=$ & $\{y\},$ &  & $\text{ }%
D_{H}(w)$ & $=$ & $H.$%
\end{tabular}%
\end{equation*}%
So, $H$ has the \emph{discrete} topology. Consequently, $H$ is $T_{6}.$
Since $H$ is finite and $T_{0},$ it follows by Lemma \ref{dim-0} (2) that $H$
is spectral. In fact, $H\approx Spec(R),$ where we may choose $R=\mathbb{F}%
\times \mathbb{K}$ for any two fields $\mathbb{F}$ and $\mathbb{K}$ (\emph{%
e.g.}, $R=\mathbb{Z}_{6}\simeq \mathbb{Z}_{2}\times \mathbb{Z}_{3}$).$%
\blacksquare $
\end{ex}

\bigskip

\begin{ex}
\label{ex-DVR}Consider the lattice $\mathcal{L}$ in Figure \ref{fig:LN}. Set
$G:=\{x,w\}$:
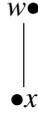
\begin{figure}[]
\centering
\begin{equation*}
\xymatrix{ w {\bullet} \ar@{-}[d] \\ {\bullet} x}
\end{equation*}%
\caption{The Sierpi\'{n}ski space is completely normal but not perfectly
normal}
\label{Sierpinski}
\end{figure}
Notice that $C^{G}(\mathcal{L})=\{x,w\}$. It follows, by Example \ref%
{normal-not-regular} and Corollary \ref{Y-top}, that $G$ attains a
Zariski-like topology induced by the poset $(L,\leq ).$ One can double
check:\ The collection of open sets of $G$ is given by%
\begin{equation*}
\begin{tabular}{lllllllllll}
$D_{G}(0)$ & $=$ & $\emptyset ,$ &  & $D_{G}(z)$ & $=$ & $\emptyset ,$ &  & $%
D_{G}(y)$ & $=$ & $\{x\},$ \\
$D_{G}(1)$ & $=$ & $G,$ &  & $D_{G}(x)$ & $=$ & $\emptyset ,$ &  & $\text{ }%
D_{G}(w)$ & $=$ & $\{x\}.$%
\end{tabular}%
\end{equation*}%
So, $G$ is the \emph{Sierpi\'{n}ski space}. Using arguments similar to the
ones above, we recover many of the well-known properties of this space, e.g.
$G$ is completely normal but not perfectly normal. Moreover, $G$ is \emph{not%
} regular (notice that $K.\dim (G)=1$). Since $X$ is finite and $T_{0},$ it
follows by Lemma \ref{dim-0} (2) that $X$ is spectral. In fact:

\begin{itemize}
\item $G\approx Spec(R),$ where $R$ is any DVR. We may choose%
\begin{equation*}
R=\mathbb{R}[[x]]\text{ with }Spec(R)=\{0,(x)\}.
\end{equation*}

\item $G\approx Spec(S),$ where $S$ is the Alaracon-Anderson semidomain%
\begin{eqnarray*}
S &=&B(n,1)\text{ with }n\geq 3,\text{ }\mathbb{P}(n-1)=\{p\}\text{ and }%
Spec(S)=\{0,pS\}; \\
S &=&B(n,n-1)\text{ with }n\geq 3\text{ and }Spec(S)=\{0,\mathfrak{m}%
_{n}\}.\blacksquare
\end{eqnarray*}
\end{itemize}
\end{ex}

\bigskip

The following table summarizes the properties of these \emph{spectral}
topological paces in the previous examples, where $\mathbf{UC:}$
ultraconnected, $\mathbf{HC:\ }$hyperconnected (irreducible), $\mathbf{R:}$
regular, $\mathbf{CR:}$ completely regular, $\mathbf{N:}$ normal, $\mathbf{%
CN:\ }$completely normal, $\mathbf{PN:\ }$perfectly normal.

{\small {%
\begin{equation*}
\begin{tabular}{|c|c|c|c|c|c|c|c|c|c|c|c|}
\hline
&  & \textbf{Graph} & $K.\dim $ & $\mathbf{UC}$ & $\mathbf{HC}$ & $\mathbf{R}
$ & $\mathbf{CR}$ & $\mathbf{N}$ & $\mathbf{CN}$ & $\mathbf{PN}$ & $T_{i}$
\\ \hline
$X$ & $\{0,x,y,w\}$ & $\mathcal{D}_{4}$ & $2$ & $\checkmark $ & $\checkmark $
& $\times $ & $\times $ & $\checkmark $ & $\times $ & $\times $ & $T_{0}$ \\
\hline
$Y$ & $\{x,y,w\}$ & $\mathcal{T}_{2}$ & $1$ & $\checkmark $ & $\times $ & $%
\times $ & $\times $ & $\checkmark $ & $\checkmark $ & $\times $ & $T_{\frac{%
3}{4}}$ \\ \hline
$Q$ & $\{x,y,w\}$ & $\mathcal{V}_{2}$ & $1$ & $\times $ & $\checkmark $ & $%
\times $ & $\times $ & $\times $ & $\times $ & $\times $ & $T_{\frac{1}{2}}$
\\ \hline
$H$ & $\{x,y\}$ & $\mathcal{P}_{2}$ & $0$ & $\checkmark $ & $\times $ & $%
\checkmark $ & $\checkmark $ & $\checkmark $ & $\checkmark $ & $\checkmark $
& $T_{6}$ \\ \hline
$G$ & $\{x,w\}$ & $\mathcal{C}_{2}$ & $1$ & $\checkmark $ & $\checkmark $ & $%
\times $ & $\times $ & $\checkmark $ & $\checkmark $ & $\times $ & $T_{\frac{%
1}{2}}$ \\ \hline
\end{tabular}%
\end{equation*}%
}}

\bigskip

We end this paper with a summary of the separation axioms for the \emph{%
Alaracon-Anderson semirings} $B(n,i)$ (cf. \ref{B(n,i)} and \ref%
{Spec(B(n,i))}). The proofs and justifications are similar to the ones in
the examples above, whence omitted.

\bigskip

\begin{exs}
\label{Table-B(n,i)}Let $n\geq 2,$ $1\leq i\leq n-1,$ $m:=n-i$ and consider
the semiring $B(n,i).$ Based on the structure of the prime spectrum as given
in \ref{Spec(B(n,i))} and applying our results in this paper, we obtain the
following summary:{\small {%
\begin{equation*}
\begin{tabular}{|c|c|c|c|c|c|c|c|c|c|c|}
\hline
$n$ & $i$ & $\omega (m)$ & $Spec(S)$ & $K.\dim (S)$ & $\mathbf{R}$ & $%
\mathbf{CR}$ & $\mathbf{N}$ & $\mathbf{CN}$ & $\mathbf{PN}$ & $T_{i}$ \\
\hline
$\geq 2$ & $0$ & $\omega (n)$ & $\mathcal{P}_{\omega (n)}$ & $0$ & $%
\checkmark $ & $\checkmark $ & $\checkmark $ & $\checkmark $ & $\checkmark $
& $T_{6}$ \\ \hline
$2$ & $1$ & $0$ & $\mathcal{P}_{1}$ & $0$ & $\checkmark $ & $\checkmark $ & $%
\checkmark $ & $\checkmark $ & $\checkmark $ & $T_{6}$ \\ \hline
$\geq 3$ & $1$ & $1$ & $\mathcal{C}_{2}$ & $1$ & $\times $ & $\times $ & $%
\checkmark $ & $\checkmark $ & $\times $ & $T_{\frac{1}{2}}$ \\ \hline
$\geq 7$ & $1$ & $\geq 2$ & $\mathcal{V}_{\omega (m)}$ & $1$ & $\times $ & $%
\times $ & $\times $ & $\times $ & $\times $ & $T_{\frac{1}{2}}$ \\ \hline
$\geq 3$ & $n-1$ & $0$ & $\mathcal{C}_{2}$ & $1$ & $\times $ & $\times $ & $%
\checkmark $ & $\checkmark $ & $\times $ & $T_{\frac{1}{2}}$ \\ \hline
$\geq 4$ & $[2,n-2]$ & $1$ & $\mathcal{C}_{3}$ & $2$ & $\times $ & $\times $
& $\checkmark $ & $\checkmark $ & $\times $ & $T_{\frac{1}{2}}$ \\ \hline
$\geq 8$ & $[2,n-2]$ & $\geq 2$ & $\mathcal{D}_{\omega (m)+2}$ & $2$ & $%
\times $ & $\times $ & $\checkmark $ & $\times $ & $\times $ & $T_{0}$ \\
\hline
\end{tabular}%
\end{equation*}%
}}
\end{exs}

\end{document}